\documentclass[12pt]{article} 
\usepackage{amsfonts,amsmath,latexsym,amssymb,mathrsfs}

\newenvironment{proof}{\noindent{\sc Proof:\ }}{\vspace{2ex}}

\def\numberlikeadb{\global\def\theequation{\thesection.\arabic{equation}}}
\numberlikeadb
\newtheorem{theorem}{Theorem}[section]
\newtheorem{lemma}[theorem]{Lemma}
\newtheorem{corollary}[theorem]{Corollary}
\newtheorem{proposition}[theorem]{Proposition}

\newcommand{\RR}{{\bf R}}

\newcommand{\PP}{{\bf P}}
\newcommand{\EE}{{\bf E}}
\def\nat{{\bf N}}
\newcommand{\var}{{\mbox{Var}}}
\newcommand{\Cov}{{\mbox{Cov}}}
\newcommand{\beas}{\begin{eqnarray*}}
\newcommand{\enas}{\end{eqnarray*}}
\newcommand{\eqs}{\begin{eqnarray*}}
\newcommand{\ens}{\end{eqnarray*}}

\newcommand{\FF}{{\cal{F}}}

\newcommand{\bea}{\begin{eqnarray}}
\newcommand{\eqa}{\begin{eqnarray}}
\newcommand{\ena}{\end{eqnarray}}
\newcommand{\eq}{\begin{equation}}
\newcommand{\en}{\end{equation}}

\def\ignore#1{}

\def\ul{^{(l)}}

\def\Ref#1{(\ref{#1})}
\def\a{\alpha}
\def\b{\beta}
\def\s{\sigma}
\def\f{\phi}

\def\l{\lambda}

\def\law{{\cal L}}

\def\ep{\hfill $\Box$ 
\def\Def{\ :=\ }

\bigskip}
\def\re{\RR}

\def\giv{\,|\,}

\def\Po{{\rm Po\,}}

\def\non{\nonumber}
\def\th{\theta}

\def\e{\varepsilon}
\def\m{\mu}

\def\var{{\rm Var\,}}

\def\r{\rho}

\def\t{\tau}
\def\g{\gamma}
\def\h{\eta}
\def\f{\phi}

\def\Bl{\left(}
\def\Br{\right)}
\def\Blm{\Bigl|}
\def\Brm{\Bigr|}
\def\Blb{\left\{}
\def\Brb{\right\}}

\def\ui{^{(1)}}

\def\Bi{{\rm Bi\,}}

\def\nin{\noindent}

\def\pr{\PP}
\def\msk{\medskip}
\def\bsk{\bigskip}

\def\ff{{\cal F}}
\def\d{\delta}
\def\z{\zeta}

\def\Eq{\ =\ }
\def\Le{\ \le\ }

\def\tH{{\widetilde H}}
\def\G{\Gamma}

\def\ignore#1{}

\def\ffo{\widetilde{\mathcal{F}}}
\def\tSk{\widetilde\Sigma^{[k]}}
\def\skk{\sum_{k=1}^K}
\def\sjj{\sum_{j=1}^J}
\def\sll{\sum_{l=1}^L}

\def\swls{\sum_{w\notin W^*_l}}

\def\uT{{}^T}
\def\uk{^{(k)}}
\def\uj{^{(j)}}

\def\uK{^{(K)}}
\def\uJ{^{(J)}}
\def\n{\nu}
\def\diag{\mbox{diag}\,}
\def\Giv{\,\Big|\,}
\def\cupdot{\cup\kern-8.2pt\cdot\kern5.5pt}
\def\cov{\Cov}
\def\tm{\tilde\m}

\def\tmm{\bar m}
\def\tnn{\bar n}

\def\tq{q}

\def\tX{\widetilde X}

\def\BA{D^X}

\def\gkxxli{G_k^{XX}(l,i)}
\def\gkxyli{G_k^{XY}(l,i)}
\def\gkxylpii{G_k^{XY}(l+1,i)}
\def\gjyxli{G_j^{YX}(l,i)}
\def\gjyyli{G_j^{YY}(l,i)}

\def\gjyxoi{G_j^{YX}(0,i)}
\def\te{{\tilde e}}
\def\tc{{\tilde c}}

\def\tth{{\tilde\th}}
\def\umn{^{(m,n)}}
\def\ggg{{\cal G}}
\def\tz{\tilde z}
\def\teps{\tilde \e}

\def\teta{\tilde \h}
\def\lb{\bar{\l}}
\def\Ekk{\EE_{k_1,k_2}}
\def\Pkk{\PP_{k_1,k_2}}
\def\bone{{\bf{1}}}
\def\erra{e(m,n)}
\def\fff{\Phi}
\def\ZZZ{\mathfrak{Z}}
\def\ZZZs{Z_*}
\def\smn{\sqrt{\frac mn}}
\def\snm{\sqrt{\frac nm}}
\def\phh{\varphi}
\def\kkk{\kappa}
\def\integ{{\mathbb Z}}
\def\tU{{\widetilde U}}
\def\tW{{\widetilde W}}
\def\nti{n\to\infty}
\def\tilt{{\tilde\t}}

\def\Def{\ :=\ }
\def\tkkk{{\tilde\kkk}}

\begin{document}

\title{The shortest distance in random multi-type intersection graphs}
\author{
A. D. Barbour\footnote{Angewandte Mathematik, Universit\"at Z\"urich,
Winterthurertrasse 190, CH-8057 Z\"URICH;  ADB was supported in part by
Schweizer Nationalfonds Projekt Nr.\ 20--117625/1, and carried out part
of the work while visiting the Institut Mittag--Leffler,
Djursholm, Sweden. 
\msk}
\ and
G. Reinert\footnote{Department of Statistics,
University of Oxford, 1 South Parks Road, OXFORD OX1 3TG, UK; GR was supported in part by EPSRC and BBSRC through OCISB.
}\\
Universit\"at Z\"urich and University of Oxford
}

\date{} 
\maketitle

\begin{abstract}
Using an associated branching process as the basis of our approximation, we show that typical inter-point distances in a multitype random intersection graph have a defective distribution, which is well described by a mixture of translated and scaled Gumbel distributions, the missing mass corresponding to the event that the vertices are not in the same component of the graph.
\end{abstract}


\section{Introduction}\label{intro}
 \setcounter{equation}{0}

\hspace{0.15in}
Bipartite graphs have been studied in a variety of applications:  directors 
and companies \cite{RobinsAlexander}, persons and questions in an intelligence 
test \cite{rasch}, or genes and gene properties \cite{tanayetal}, to give 
just a few examples. Typically, in such applications, the graph induced on
the vertices of one of the two parts, with vertices linked if there is a path
of length~2 joining them in the bipartite graph, are of primary interest.
For instance, the structure of the network linking directors may be of
greater interest than the companies involved.  Furthermore, in some applications,
the remaining part of the bipartite graph, which is responsible for
forming the links, may not be known or observable, and it may be of interest
to deduce its existence from the properties of the observed part of the
structure alone. However, the statistical properties of such bipartite graphs are 
not well understood, particularly when there are different types of vertices, 
see \cite{RobinsAlexander}. Here, we shall be concerned with the properties of a
particular family of such graphs, known as random intersection graphs, and
with the statistics of distances between randomly chosen points. 

Random intersection graphs are constructed from two sets, the `vertices'
and the `objects', as follows.
Each vertex~$v \in V$  is associated with a randomly chosen subset~$A_v$
of a finite set~$U$ of objects, and two vertices $v$ and~$v'$ are joined in 
the graph if $A_v \cap A_{v'} \neq \emptyset$.   In the simplest case,
the Bernoulli model, 
vertex~$v$ is associated with object~$u$ independently of all other associations
with fixed probability~$p$.  
Britton {\it et al.} \cite{BDLL} establish a branching process approximation 
for the spread of a Reed-Frost epidemic on such a graph.  Here, we consider the more flexible
model in which there are~$K$ distinct types of vertices and~$J$ types of
objects, and vertex~$v$ of type~$k$ is associated with object~$u$ of type~$j$
independently of all other associations with probability~$p_{kj}$. Our model 
can be viewed as a bipartite Erd\H os--R\'enyi mixture graph \cite{daudinetal}.
In Erd\H os--R\'enyi mixture graphs, vertices are coloured, with
the probability of two vertices being connected depending only on their colours;
edges occur independently. 
 
Random intersection graphs of this kind can also be related to the Rasch \cite{rasch}
 models in social science.  These are given by taking 
$$
   p_{kj} \Eq \frac{\alpha_k \beta_j}{1 + \alpha_k \beta_j}.
$$
For example, one might have $k \sim j$ if person $k$ solves problem~$j$ 
correctly; the $\alpha$'s would then relate to the ability of the person, 
and the $\beta$'s to the type of problem. A simplified Rasch model of the form 
$   p_{kj} \Eq {\alpha_k \beta_j}$
can be viewed as a special case of an exponential random graph model, see Equation (1) in \cite{robinsetal}.


In the study of random networks, the shortest distance between two randomly 
chosen vertices is one of the standard summary statistics. In this paper, 
we approximate its distribution for multitype Bernoulli intersection graphs. 
Since the networks used in applications are typically finite, we not only 
provide a limiting approximation, but also give explicit bounds on the 
difference between the true and limiting distributions.  Our main results,
summarized in Corollary~\ref{Main}, give an approximation described
in terms of the distributions of the limiting random variables~$W$ of the
associated multivariate bipartite branching process, when the process starts
with a single individual of one or other of the types.  The probability of
the two vertices being in the same component is well approximated by the product
of the probabilities that neither of the branching processes becomes extinct.
On this event, the distance has a distribution close to that of a 
translation--mixture of scaled Gumbel distributions, with the mixture
distribution being explicitly given in terms of those of the limiting
random variables~$W$.  Alternatively, the approximate distribution
can be described as that of (a linear transformation of) the sum of 
three independent random 
variables, one a Gumbel, and the others distributed as the logarithm of~$W$,
given the appropriate initial types. In  a natural asymptotic
framework, the error bounds behave like an inverse power of the total
number of vertices, whose exponent can be derived from the parameters
of the bipartite graph: the probabilities $p_{kj}$, and the numbers $n_k$ of
vertices, $1\le k\le K$, and $m_j$ of objects, $1\le j\le J$, of the
different types.

The structure of the paper is as follows.  The link between intersection
graphs and branching processes is described in Section~\ref{IB}.
The necessary distributional properties of the branching process are
established in Section~\ref{offspring}, and the extent to which it
differs from the intersection graph is controlled in Section~\ref{Ghosts}.
The main theorem is then stated and proved in Section~\ref{distances},
and an application to exponential random graph models is given in Section~\ref{ER}.  A key element in the
proof is a Poisson approximation to coincidence probabilities in a
generalization of the hypergeometric sampling scheme; this is
undertaken in Section~\ref{intersection}.

Although our motivation for studying this problem comes from the
bipartite setting, one could equally well conduct a similar analysis
for a graph without bipartite structure, recovering the general
Erd\H os--R\'enyi mixture model; a corresponding approximation
is given without detailed proof in~\Ref{ER-result}.  
However, the analysis for a `general'
graph would not easily imply our results as a special case,
with the vertices split into two groups and with a bipartite
matrix of edge probabilities~$P$, because the 2-periodic structure would result
in there being more than one eigenvalue of the mean matrix
having largest modulus, and methods such as those of this paper would still be
needed, to deal with the extra complexity that results.

\section{Intersection graphs and branching processes}\label{IB}
 \setcounter{equation}{0}
A random multitype intersection graph on the vertex set 
$V = V_1\cupdot\ldots\cupdot V_K$ is defined using a second set of
`objects'
$U = U_1\cupdot\ldots\cupdot U_J$: each vertex $v \in V_k$
independently chooses a subset $A_v \in U$ with distribution depending 
on~$k$ alone, and $v \sim v'$ if and only if $A_v \cap A_{v'} \neq \emptyset$.  
Here, we restrict ourselves to graphs derived from an underlying
Erd\H os--R\'enyi bipartite mixture model, in which only 
edges~$e_{uv}$ between $u\in U$ and $v\in V$ are possible, and these are 
present or absent independently, with probability $p_{kj}$ if $u\in U_j$ 
and $v\in V_k$.  Thus, in the random intersection graph itself, $v\sim v'$ if 
and only if, for some $u\in U$, both $e_{uv}$ and $e_{uv'}$ are present.   
 
Such a random graph can be constructed from a bipartite multitype
branching process $(Z(0),Z(1),Z(2),\ldots) = (X(0),Y(1),X(1),\ldots)$,
with $X(i) \in \nat^K$ and $Y(i) \in \nat^J$ for each~$i$, together
with sets of randomly assigned indices.  Start with numbers
$X(0) = (X_1(0),\ldots,X_K(0))$ of individuals of the different types
$\{(k,1),\,1\le k\le K\}$. 
The \hbox{$s$-th} individual (in some ordering) of type~$(k,1)$ in $Z(2r)
= X(r)$ has offspring vector $Y_{ks;r} = (Y_{ks;r}\ui,\ldots,Y_{ks;r}\uJ)$,
realized from the product distribution $\otimes_{j=1}^J \Bi(m_j,p_{kj})$,
where $m_j$ is the cardinality of the set~$U_j$, and the random vectors
$(Y_{ks;r},\,1\le k\le K,s\ge1,r\ge0)$ are independent.  We then set
\eq\label{Yksr-def}
    Z(2r+1) \Eq Y(r+1) \Eq \skk  \sum_{s=1}^{X_k(r)} Y_{ks;r}.
\en
Similarly, the $t$-th individual of type $(j,2)$ in $Z(2r-1) = Y(r)$ has
offspring vector $X_{jt;r} = (X_{jt;r}\ui,\ldots,X_{jt;r}\uK)$,	
realized from the product distribution $\otimes_{k=1}^K \Bi(n_k,p_{kj})$,
where $n_k$ is the cardinality of the set~$V_k$, and the random vectors
$(X_{jt;r},\,1\le j\le J,s\ge1,r\ge1)$ are independent of each other 
and of the $Y_{ks;r}$.	We then set
\[
    Z(2r) \Eq X(r) \Eq \sjj  \sum_{t=1}^{Y_j(r)} X_{jt;r}.
\]
We also define 
\eq\label{MN-def}
   m \ :=\ \sjj m_j; \qquad n \  :=\ \skk n_k.
\en
Throughout this paper we assume that $m_j \ge 2$ for $ 1 \le j \le J$, and $n_k \ge 2$ for $1 \le k \le K$, and that $X_k(0) \le n_k$ for $1 \le k \le K$. 

To obtain the intersection graph, label each individual in the bipartite
branching process with its line of descent.  Thus
\[
   \{i; (k_0,s_0),(j_1,t_1),(k_1,s_1),\ldots,(j_i,t_i),(k_i,s_i)\}
\]
labels the $s_i$-th individual of type $(k_i,1)$ in generation~$2i$,
which was descended from the $t_i$-th individual of type $(j_i,2)$
in generation~$2i-1$, and so on.  These labels are then augmented
with indices from the index set appropriate to the type of individual, 
as follows. The $Y_{k_i s_i;i}\uj$ type $(j,2)$ offspring of the typical
vertex above are each assigned at random a unique index from a
uniformly and independently chosen subset $\law_{k_i s_i;i}\uj \subset
\{1,2,\ldots,m_j\}$ of size $Y_{k_i s_i;i}\uj$; a similar construction
is used for the offspring of objects.  

  A further class identifier,
$0$ or~$1$, is then attached to each individual: $1$ if the individual
is in generation zero, and thereafter, taking the individuals of the 
bipartite process in order of generation, but in any order within each
generation, assign class~$0$ if its parent was in class~$0$,
or if its index had previously been assigned to another individual of
the same type and of class~$1$, and~$1$ otherwise.  
The class~$0$ individuals we refer to as ghosts. Edges are also
created between the indices of a parent and its child if both are 
of class~$1$, or
if the parent is of class~$1$ and the child of class~$0'$, 
where class~$0'$ indicates a class~$0$ individual whose
index was first assigned (therefore to a class~$1$ individual)
in its own generation. Then the individuals that belong
to class~$1$ correspond, via their indices, to the vertices
and objects used in constructing the intersection graph, and
two vertices have an edge between them if there are corresponding 
class~$1$ or class~$0'$
individuals in the bipartite branching process that are at distance~$2$
from one another.  In this way, the union of the components of the 
random intersection graph that contain the initial vertices is
sequentially constructed according to distance from the initial
vertices, the
class~$1$ vertices in generation~$2i$ of the bipartite branching 
process corresponding to the vertices that are at distance~$i$ in the 
random intersection graph from the initially chosen set of vertices.
If these components do not exhaust all vertices, the process can
be continued from any unused vertex, until all have been used.

Two vertices,  $A$ of type $k_1$ and $B$ of type~$k_2$,
are at distance at least~$d+1$ from one another in the random
intersection graph if the $d$-neighbourhood
of one of them in the bipartite graph does not intersect the 
$d$-neighbourhood of the other.  
Constructing the random intersection graph from a bipartite branching 
process starting with one individual~$A$ of type $(k_1,1)$ and one, $B$, of type
$(k_2,1)$, this is the case exactly on the event that the set of all 
class~$1$ or class~$0'$ descendants of~$A$ up to time~$d$ --- both of 
types $(k,1)$ and of types $(j,2)$ ---  is disjoint from that of~$B$.
From the construction, these
sets can only overlap if there is at least one class~$1$
descendant of either $A$ or~$B$ having the same index as a class~$0'$
descendant of the other, and then necessarily in the same generation
of the bipartite process.
Our main theorem consists of showing that
the probability of this event can be well approximated by the
probability of the corresponding event when {\it all\/} descendants
are considered, and that this probability in turn can be well approximated using 
the theory of branching processes. 
 
The origin of the approximation lies in the following well known facts 
(\cite{Harris}, II Theorem 9.2):
that,  on the event of non-extinction, a square integrable super-critical 
multitype branching process, whose mean matrix is irreducible and aperiodic,
has an asymptotically stable type distribution; 
and that the total number of individuals alive in each generation grows
like a random multiple of a geometrically growing sequence.
For the~$X$ branching process, this means that the number of individuals
of type~$k$ in generation~$i$ is approximately given by $\t^i W \m_k$,
where~$\t$ is the largest eigenvalue of the mean matrix, $\m^T$ is the
associated positive left eigenvector, and~$W$ is a non-negative random 
variable, positive on the event of non-extinction, and the same for all
$i$ and~$k$.  Hence the numbers of descendants~$X^A(i)$ of~$A$ at the $i$-th 
generation of the~$X$ branching process are approximately given by $\t^i W^A \m_k$,
and those of~$B$ by $\t^i W^B \m_k$, where $W^A$ and~$W^B$ are independent.
When constructing the random intersection graph from the branching process,
indices are assigned to the vertices independently at random, with replacement.
Links between the $A$ and~$B$ neighbourhoods occur whenever, for some $i\ge1$
and some $1\le k\le K$, one or more of the~$X_k^B(i)$ are assigned the same
index as one of the~$X^A(i)$; other coincidences give rise to `ghosts',
and play no part in the intersection graph.  The mean number of such events 
up to and including generation~$i$ is thus approximately
\[
   \sum_{s=1}^i \t^{2s} W^A W^B \skk \frac{\m_k^2}{n_k} \ \approx\ 
    \kkk^X n^{-1}\t^{2i} W^AW^B,
\]
where
\[
     \kkk^X \Def \Bl\frac{\t^{2}}{\t^2-1}\Br \skk \frac{\m_k^2}{q^X_k},
\]
and $q^X_k := n_k/n$.
A similar formula hold for links occurring because of coincidence of indices
at the object level; here, the expected number of links up to and including
generation $i-1$ is approximately $\kkk^Y n^{-1}\t^{2(i-1)} W^AW^B$, where, 
because of~\Ref{z-reln}, $\kkk^Y = \t\kkk^X$.  Adding the two means gives
an overall mean number of links approximately equal to $\kkk n^{-1}\t^{2i} W^AW^B$,
where 
\eq\label{kappa-def}
  \kkk \Def \kkk^X(1 + \t^{-1}) \Eq \frac\t{\t-1}  \skk \frac{\m_k^2}{n_k}.
\en
Then, using Poisson
approximation, it follows that the probability of there being no shared
vertices in the $i$-neighbourhoods of $A$ and~$B$ is approximately
\[
   \Ekk\Blb e^{-\kkk n^{-1}\t^{2i} W^AW^B} \Brb,
\] 
this being the probability that the distance between $A$ and~$B$ 
in the intersection graph exceeds~$2i$.
This line of reasoning is made precise in the coming sections, and the detailed
results are to be found in Theorem~\ref{ip-distance} and Corollary~\ref{Main}.

\section{Counting the offspring}\label{offspring}
\setcounter{equation}{0}
We now study the bipartite branching process~$Z$ in greater detail. 
Our aim in this section is to justify the simple approximation,
outlined above, to the numbers $X_k(i)$ of type-$(k,1)$ individuals 
in~$Z(2i)$ (or, equivalently, of type-$k$ vertices in the $i$-th 
generation of the vertex branching process) and~$Y_j(i)$ of 
type-$(j,2)$ individuals in $Z(2i-1)$,
(or of objects in the $i$-th generation of the object branching process). 
Theorems \ref{TX-W-bnd} and~\ref{TY-W-bnd} below show that, for large~$i$,
$X(i) \sim \t^i W\m$ and $Y(i) \sim \t^{i-1}\z W\tm$, 
the notation being as defined below.  

\subsection{Assumptions and notation}\label{notation}
Let
\beas
   N_X &=& \diag(n_1, \ldots, n_{K}) \quad \mbox{and} 
            \quad N_Y \Eq \diag(m_1, \ldots, m_{J}) 
\enas
be the diagonal matrices of the numbers of different types of vertices, 
and of different types of objects, respectively; for future convenience,
recalling~\Ref{MN-def}, we define
\eq\label{q-defs}
   q_k^X\ :=\ \frac{n_k}{n},\quad 1\le k\le K;\qquad
   q_j^Y\ :=\ \frac{m_j}{m},\quad 1\le j\le J.
\en
Let $P$ denote the $K \times J$-matrix of edge probabilities, and put
$$
   M_X \Eq  P N_Y P^T N_X.
$$
Then the non-negative matrix~$M_X$ is the mean matrix of the~$X$ branching
process.

\medskip
\nin {\bf Assumption.}\quad We assume that the non-negative matrix~$M_X$ is 
irreducible and aperiodic, and has largest eigenvalue~$\t > 1$. 

\medskip
We use $\n$ and~$\m^T$ to denote respectively the right and left eigenvectors corresponding
to~$\t$,  with $\mu_k > 0, \, 1 \le k \le K$, 
standardized so that $\|\mu\|_1 =1$ and that $\mu^T \nu = 1$. 
We assume throughout that $\t>1$. We then define 
$M_Y :=  P^T N_X P N_Y$ to be the mean matrix of the~$Y$ branching process,
and
\eq\label{tm-defs} \tm :=  \z^{-1} N_YP^T\m
\en  to be the left eigenvector of~$M_Y$ with 
$\|\tm\|_1 = 1$
corresponding to the eigenvalue~$\t$. Thus $\z := \mu^TPN_Y \bone$, where $\bone$ 
denotes a $J \times 1$-vector of 1's.

Note also that 
\beas
 \z^2 \sjj \frac{\tm_j^2}{m_j} &=& (N_Y P^T \m)^T N_Y^{-1} (N_Y P^T \m) 
      \Eq \mu^T P N_Y P^T \m\\
 &=&   \mu^T M_X N_X^{-1} \m  \Eq \tau \mu^T  N_X^{-1} \m 
      \Eq \t \skk \frac{\m_k^2}{n_k} , 
\enas 
so that 
\eq\label{z-reln}
  \z^2 n/m \Eq \t \Blb \skk  \frac{\m_k^2}{q_k^X} \Brb
         \Big/ \Blb \sjj \frac{\tm_j^2}{q_j^Y} \Brb \ =:\ \ZZZ^2,
\en
say.

We next define~$c_0$ to be the smallest value such that
\eqa
  \sup_{a\colon \|a\|_1 = 1} \sup_{i\ge0} \t^{-i}a^T M_X^i e\uk 
	   \Le c_0\m_k, \qquad 1\le k\le K; \non\\
  \sup_{a\colon \|a\|_1 = 1} \sup_{i\ge0} \t^{-i}a^T M_Y^i \te\uj 
	   \Le c_0\tm_j, \qquad 1\le j\le J, \label{A}
\ena
where $e\uk$ and~$\te\uj$ denote the $k$ and~$j$ unit vectors 
in $\re^K$ and~$\re^J$ respectively.
Further, with $\lambda_2$ the eigenvalue of $M_X$ with second largest modulus, 
we define~$c_1$ such that
\eq\label{A-2}
  \sup_{b\colon \|b\|_\infty = 1} \sup_{i\ge0} 
	       |\l_2|^{-i}\|(M_X - \tau \nu \mu^T)^i b\|_\infty \Le c_1, 
\en
where we take $(M_X - \tau \nu \mu^T)^0 := I - \n\m^T$. 
Note that it follows from the Perron-Frobenius Theorem that  $c_0$ and~$c_1$ are both finite; 
see \cite{HornJohnson}, Theorem~8.5.1, and \cite{mode}, Chapter~1, Theorem~6.1.  We also, 
for later use, write
\eq\label{gamma-eta} 
    \g\ :=\ \max\{\t,|\l_2|^2\} < \t^2;\qquad 
    \th\ :=\ \max_{s\ge0} (s+1)(\sqrt\g/\t)^s,
\en	
and introduce the notation $\ffo_r$ to denote the $\s$-algebra $\s\{Z(t),\,0\le t\le r\}$,
and $\FF_i^X:= \s\{X(l),\,0\le l\le i\}$, $\FF_i^Y:= \s\{Y(l),\,1\le l\le i\}$.

\subsection{Asymptotics}\label{asymptotics}
The main results of the paper require no particular asymptotic setting. 
However, asymptotics are useful for putting the results in the context
of a natural limiting framework.  One such
choice is the following.  Start by choosing the $n_k$ 
and~$m_j$ so that the proportions $\tq_k^X(m,n)$ and~$\tq_j^Y(m,n)$ converge to
non-zero limits. Then one can arrange for~$P = P\umn$ to vary as $m$ and~$n$ tend to 
infinity, in such a way that the matrix~$M_X\umn$ converges to a fixed irreducible and
aperiodic matrix~$M_X$, entailing the convergence of quantities such
as $\t\umn$, $\m\umn$ and~$\n\umn$ to limits $\m$, $\n$ and~$\t$.  
With this in mind, define $Q_X\umn := n^{-1}N_X$ and 
$Q_Y\umn := m^{-1}N_Y$, and take $P\umn := (mn)^{-1/2}\Pi$, for a fixed matrix~$\Pi$.
This then gives $M_X\umn = \Pi Q_Y\umn \Pi^T Q_X\umn$,
so that, if $Q_X\umn \to Q_X$ and $Q_Y\umn \to Q_Y$, with $Q_X$ and~$Q_Y$
both having positive diagonals, then $M_X\umn \to M_X := \Pi Q_Y \Pi^T Q_X$.
If also, in keeping with the general assumptions of the paper, we have $\t > 1$,
then we describe this behaviour as `standard asymptotics'.  

  Other asymptotic settings could
equally well be considered.  For instance, there would be no great difference
in the qualitative behaviour if~$\t\umn$ were allowed to tend to infinity 
with~$n$ like a power of~$\log n$.

\subsection{Expectations}\label{expectations}
We begin our analysis by examining the growth of the mean
numbers of individuals of different types.
Using $\EE_0$ to denote $\EE\{\cdot\giv \ffo_0\}$, we 
immediately have
\eq\label{X-mean}
   \EE \{X^T(i) \giv \ffo_{2i-1}\} \Eq Y^T(i) P^T N_X;
	       \qquad \EE_0 X^T(i) \Eq X^T(0) M_X^i,
\en
and 
\eq\label{Y-mean} 
   \EE \{Y^T(i) \giv \ffo_{2i-2}\} \Eq   X^T(i-1)P N_Y;\qquad  
	                          \EE_0 Y^T(i) \Eq X^T(0) M_X^{i-1}PN_Y .
\en
From these, and using~\Ref{A}, we have, for instance,
\beas
    \EE\{X_k(i) \giv \FF_s^Y\} 
      &=& \sjj Y^T(s)M_Y^{i-s}\te\uj\,(\te\uj)^T P^TN_Xe\uk \\
  &\le& c_0\|Y(s)\|_1 \t^{i-s}\tm^TP^TN_Xe\uk 
     \Eq  c_0\|Y(s)\|_1 \t^{i-s}\z^{-1}\m^T M_X e\uk \\
  &=& c_0\|Y(s)\|_1 \t^{i-s}(\t/\z)\m_k,
\enas
so that, for any $1\le k\le K$ and $1\le j\le J$, and for $i \ge s \ge 0$,
\bea
  \EE\{X_k(i) \giv \FF_s^X\}  \Le c_0 \|X(s)\|_1 \t^{i-s} \m_k;&& \hspace{-4.5mm} 
     \EE\{X_k(i) \giv \FF_s^Y\} \Le c_0 \|Y(s)\|_1 \z^{-1}\t^{i-s+1} \m_k;\phantom{H}
				\label{A0} \\
 \EE\{Y_j(i) \giv \FF_s^Y\}  \Le c_0 \|Y(s)\|_1 \t^{i-s} \tm_j;				
   &&  \hspace{-4.5mm} 
	  \EE\{Y_j(i) \giv \FF_{s-1}^X\} \Le c_0 \|X(s-1)\|_1 \z\t^{i-s} \tm_j;
				\label{A01} \\			
    \EE_0 X_k(i) \Le c_0 \|X(0)\|_1 \t^i \m_k;&& \hspace{-4.5 mm} 
	      \EE_0 Y_j(i) \Le c_0 \|X(0)\|_1 \z \t^{i-1} \tm_j. \label{A1}
\ena

\subsection{X-Covariances}\label{X-covariances}
Controlling the covariances of the components of~$X(i)$ and~$Y(i)$ requires more work.
To start with, we observe that 
\beas
   \EE (X^T(i+1) \n \giv \FF_i^X) &=& X^T(i) M \n \Eq \t X^T(i) \n,
\enas 
so that $W_i := \t^{-i} X^T(i) \n$, $i\ge0$, is a non-negative martingale with respect to the filtration~$\{\FF_i^X, i=0,1, \ldots \}$
which converges almost surely to a limit~$W$, and 
\bea \label{1} 
   \EE_0 ( X^T(i)\n) \Eq \t^{i} X^T(0) \n. 
\ena
The variability in the branching process is essentially determined by that 
of~$W$, which is itself largely determined during the early stages of
development.  Indeed, writing $X^T(i) = X^T(i)\{I-\n\m^T\} + X^T(i)\n\m^T$, 
we show that the variance of $X^T(i)b$ is dominated by 
$(\m^T b)^2 \var(X^T(i)\n)$, unless $\m^Tb = 0$.

\begin{lemma}\label{mart-varce}
The variance of the martingale $W_i$ is bounded as follows:
\bea
   (i)&& \var_0 \left( W_i - W \right) 
		        \Le c_2 \{\t/(\t-1)\} \|X(0)\|_1 \t^{-i}; \phantom{HHHHHHH} \non\\ 
   (ii)&& \var_0 W_i \Le c_2\{\t/(\t-1)\} \|X(0)\|_1 ,
	      \non 
\ena 
where $c_2 := c_0 \t^{-2}\|\n\|_\infty^2 S(\Sigma)$ and  $S(\Sigma) := \max_{1\le k\le K} \sum_{l,m=1}^K|\Sigma^{[k]}_{lm}|$: $\Sigma^{[k]}$ is defined below.
\end{lemma}

\proof
We begin by writing
\beas
     X^T(i+1) &=& \sum_{k=1}^K \sum_{r=1}^{X_{k}(i)}  \tX^T_{k r; i} ,
\enas
where $\tX_{k r; i}$ denotes the $K$-vector of descendants in $X$-generation $i+1$
of individual~$r$ of type~$k$ in $X$-generation~$i$.  
The random vectors
$(\tX_{k r; i},\,1\le k\le K,r\ge1,i\ge0 )$ are independent, and, 
for each $k$, the $\tX_{k r; i}$ are identically distributed, with means the 
transpose $M_{X,[k]}$ of the $k$-th row of $M_X$, and with a covariance matrix
that we denote by~$\Sigma^{[k]}$. Then
\beas
   X^T(i+1) &=& \sum_{k=1}^K \sum_{r=1}^{X_{k}(i)} M^T_{X,[k]} 
	    + \sum_{k=1}^K \sum_{r=1}^{X_{k}(i)} ( \tX^T_{k r; i} - M^T_{X,[k]} ),
\enas
and so
\eq \label{ADB-1}
   X^T(i+1) - X^T(i) M_X 
     \Eq \sum_{k=1}^K \sum_{r=1}^{X_{k}(i)} (\tX^T_{k r; i} - M^T_{X,[k]}). 
\en 
Considering the right hand side, we have
\bea
  &&\EE \left\{\left(\sum_{k=1}^K \sum_{r=1}^{X_{k}(i)} (\tX_{k r; i} - M_{X,[k]})\right)
    \left(\sum_{k'=1}^K \sum_{r'=1}^{X_{k'}(i)} (\tX_{k' r'; i} - M_{X,[k']})\right)^T 
	       \,\bigg|\, \FF_i  \right\}  \non\\
   &&\Eq \sum_{k=1}^K \sum_{r=1}^{X_{k}(i)} 
	\EE\{(\tX_{k r; i} - M_{X,[k]})(\tX_{k r; i} - M_{X,[k]})^T \giv \FF_i^X\} \non\\
   &&\Eq  \sum_{k=1}^K X_k(i) \Sigma^{[k]}, \label{star}
\ena
by the independence of the vectors and their having mean zero.
Hence, using (\ref{star}), it follows that
\bea  
    \EE \left\{ \bigl(\t^{-(i+1) }X^T(i+1)\n - \t^{-i}X^T(i)\n \bigr)^2  
	           \Giv \FF_i \right\} 
    &=& \t^{-2(i+1)} \sum_{k=1}^K X_k(i) \n^T\Sigma^{[k]}\n , \label{4}
\ena
and thus, from~\Ref{A1},
\bea 
   &&\EE_0 \left\{ \bigl(\t^{-(i+1)}X^T(i+1)\n - \t^{-i}X^T(i)\n \bigr)^2 \right\} \Le 
	    \sum_{k=1}^K c_0 \|X(0)\|_1 \t^{-i-2} \m_k \n^T\Sigma^{[k]}\n \phantom{HHH}\non\\
   &&\qquad \Le c_0 \|X(0)\|_1 \t^{-i-2} \|\n\|_\infty^2 S(\Sigma) 
        \Eq c_2 \|X(0)\|_1 \t^{-i}, \label{var-1}
\ena
where $S(\Sigma) := \max_{1\le k\le K} \sum_{l,m=1}^K|\Sigma^{[k]}_{lm}|$
and $c_2 := c_0 \t^{-2}\|\n\|_\infty^2 S(\Sigma)$.
Writing the martingale  $\t^{-i} X^T(i)\n$ as a sum of its one-step differences,
the lemma now follows easily. \ep

The next lemma controls the variances of those components of~$X(i)$ that are 
orthogonal to~$\n$; note that
\[
     X^T(i) \Eq (X^T(i)\n)\m^T + X^T(i)(I - \n\m^T).
\]  

\begin{lemma}\label{orthog-varce}
With $c_2$ as in Lemma~\ref{mart-varce} and $\g := \max\{\t,|\l_2|^2\} < \t^2$, and with $c_3 := c_1^2 S(\Sigma)$, 
we have
\bea
   (i)&& \var_0\{X^T(i)(I-\n\m^T)b\} \Le c_3 \|X(0)\|_1 \,i\g^i\, \|b\|_\infty^2; 
     \non\\ 
  (ii)&& \var_0\{X^T(i)b\} \Le 2\|X(0)\|_1\,
       \{c_2 \{\t/(\t-1)\}  \,\t^{2i}\, (\m^Tb)^2 +  c_3 \, i\g^i\, \|b\|_\infty^2 \} \phantom{HH}
     \non 
\ena
for any $b\in\re^K$.  In particular, with~$b = e_k$ the~$k^{th}$ unit vector it follows  that
\[  
 (iii)\qquad  \var_0\{X_k(i) - (X^T(i)\n)\m_k\} \Le c_3 \|X(0)\|_1\, i\g^i,\qquad 1\le k\le K.
   \phantom{HHHHH}
\]  
\end{lemma}

\proof
Recalling that $\mu^T \n = 1$, it can be seen by
induction that
\bea
    \left(M_X - \t (\n \mu^T) \right)^i &=& M_X^i - \t^{i} (\n \mu^T). \label{p7star} 
\ena 
Hence, with~$M_X^0 = I$ and as~$\n \mu^T (M_X - \tau \n \mu^T)=0$, 
\bea
   \lefteqn{ X^T(i) ( I - \n \mu^T) } \non\\
   &=& \sum_{r=0}^{i-1} \left\{ X^T(r+1) \left[M_X^{i - r - 1} - (\n\mu^T) \t^{(i-r-1)} 
      \right]
	     - X^T(r) \left[M_X^{i - r} - (\n \mu^T) \t^{(i-r)} \right]  \right\}\non \\
   &&\mbox{}+ X^T(0) \left[M_X^{i} - (\n \mu^T) \t^{i} \right]  \non\\
	 &=&  \sum_{r=0}^{i-1} \left\{ X^T(r+1)(I - \n \mu^T) 
           - X^T(r) (M_X - \t (\n \mu^T)) \right\} 
       \left[M_X - \t (\n \mu^T)\right]^{i-r-1} \non\\
   && + X^T(0) \left[M_X^{i} - (\n \mu^T) \t^{i} \right] \non\\
   &=&  \sum_{r=0}^{i-1}  U^T(r) A^{i-r-1} 
      + X^T(0) \left[M_X^{i} - (\n \mu^T) \t^{i} \right], \label{var-2}
\ena 
where 
\[
    U^T(r) \ :=\ (X^T(r+1) - X^T(r)M_X) (I - \n \mu^T) \quad \mbox{and}
		  \quad A \ :=\ M_X - \t (\n \mu^T).
\]
Note, in particular, that $\EE_0 \{U(r)U^T(s)\} = 0$ whenever $r\neq s$, in view
of~\Ref{ADB-1}, and that
\eq \label{var-3}
    \EE \{U(r)U^T(r) \giv \FF_r^X\} \Eq \sum_{k=1}^K X_k(r)(I-\m\n^T)\Sigma^{[k]}(I-\n\m^T).
\en		
Hence, since $(I-\n\m^T)(M_X-\t\n\m^T) = M_X-\t\n\m^T$, it follows
from~\Ref{A-2} and~\Ref{A1} that
\bea
  \lefteqn{\var_0\{X^T(i)(I-\n\m^T)b\} \Eq b^T\sum_{r=0}^{i-1} (A^T)^{i-1-r} 
        \EE_0\{U(r)U^T(r)\}   A^{i-1-r}b } \non\\
  &&\Le    \sum_{r=0}^{i-1} \sum_{k=1}^K \EE_0 X_k(r)\{c_1 |\l_2|^{(i-1-r)}
	   \|b\|_\infty\}^2 S(\Sigma) \phantom{HHHHH}\non\\
  &&\Le c_0\|X(0)\|_1 \sum_{r=0}^{i-1} \t^r \{c_1 |\l_2|^{(i-1-r)}
	   \|b\|_\infty\}^2 S(\Sigma), 
\ena
for any $b \in \re^K$, proving part~(i), and part~(iii) follows directly.
Since also, from Lemma~\ref{mart-varce}\,(ii),   
\eq\label{var-4a}
   \var_0\{X^T(i)\n\m^Tb\} \Le c_2 \{\t/(\t-1)\} \|X(0)\|_1 \t^{2i} (\m^Tb)^2,
\en
it follows from part~(i) that
\bea
   \var_0\{X^T(i)b\} &\le&  2\bigl\{\var_0\{X^T(i)\n\m^Tb\} 
                             + \var_0\{X^T(i)(I-\n\m^T)b\}\bigr\} \phantom{HHHHHHH}\non\\
      &\le& 2\{c_2 \{\t/(\t-1)\} \|X(0)\|_1 \,\t^{2i}\, (\m^Tb)^2 + 
                       c_3 \|X(0)\|_1\, i\g^i\, \|b\|_\infty^2 \}, 
\ena
establishing part~(ii).  \ep

\nin Note that the growth of $\var_0\{X^T(i)b\}$ with~$i$ is at
 rate $O(i\g^i\, \|b\|_\infty^2)$, slower than~$\tau^{2i}$,  if $\m^Tb = 0$. Note also
that, if $\t \neq |\l_2|^2$, the factor~$i$ can be replaced by a constant
$c(\t,|\l_2|^2)$.

\begin{corollary}\label{product-moment}	 	  
For all $1 \le k,l \le K$, there are constants $c_5$, $c_6$ such that
\bea
  \EE_0 \{X_k(i)X_l(i)\} 
	  &\le& c_5 \|X(0)\|_1^2 \m_k\m_l \t^{2i} + c_6 \|X(0)\|_1 i \t^i\g^{i/2}.\non
\ena
In particular, for $c_5' := c_5 + c_6K^2\th$, where~$\th$ is as in~\Ref{gamma-eta},
we have
\[
  \EE_0 \|X(i)\|_1^2 \Le c_5' \t^{2i} \|X(0)\|_1^2.
\]
\end{corollary}

\proof
It follows from \Ref{A1} and Lemma~\ref{orthog-varce}\,(ii) with the Cauchy-Schwarz inequality that
\bea
  \lefteqn{\EE_0 \{X_k(i)X_l(i)\} \Le |\EE_0 X_k(i) \EE_0 X_l(i)| + |\cov_0(X_k(i),X_l(i))|} \non\\ 
  &\le&   \{c_0^2 \|X(0)\|_1^2  + 2c_2\|X(0)\|_1 \t/(\t-1)\} \m_k\m_l \t^{2i} \non\\
  &&\qquad\mbox{}+ 	
      4\|X(0)\|_1 (i\g^i)^{1/2}\t^i \sqrt{c_2 c_3 \{\t/(\t-1)\}} + 2c_3 \|X(0)\|_1 i\g^i,
\ena
and the corollary follows by taking $c_5 := c_0^2 + 2c_2\t/(\t-1)$ and
$c_6 := 4\sqrt{c_2 c_3 \{\t/(\t-1)\}} + 2c_3$.   \ep

\subsection{Y-Covariances}\label{Y-covariances}
Very similar arguments can also be carried through for the vectors~$Y(i)$, $i\ge1$.
We first show that $Y_j(i)$ is close enough to $X^T(i-1)PN_Y \te\uj$, for any~$1\le j\le J$.
Let $Z_j(i) := Y_j(i) - X^T(i-1)PN_Y \te\uj$.

\begin{lemma}\label{Y-bounds}
There is a constant~$c_{10}$ such that, for $1\le j\le J$,
\eq
  \var_0\{Z_j(i)\} \Le c_{10}\|X(0)\|_1 \t^{i}\z\tm_j. \label{var-Y-compts}
\en
\end{lemma}

\proof
First, the quantity $Z^T(i) := Y^T(i) - X^T(i-1)PN_Y$ is represented in fashion analogous 
to~\Ref{ADB-1}. For independent random vectors $Y_{ks, r}$ as in~\Ref{Yksr-def}, 
\[
Z^T(i) =  Y^T(i) - X^T(i-1)PN_Y = \sum_{k=1}^K \sum_{s=1}^{X_k(i-1)} (Y_{ks, i-1} - \EE Y_{ks, i-1}), 
\] 
 from which it follows that
\[
   \EE\{Z(i)Z^T(i) \giv \FF_{i-1}^X\}
    \Eq \skk X_k(i-1) \tSk,
\]
where $\tSk := \diag\{m_jp_{kj}(1-p_{kj}),\,1\le j\le J\}$, and hence that
\[
   \EE_0\{Z(i)Z^T(i)\} \Eq \skk \EE_0 X_k(i-1) \tSk.
\]
This with~\Ref{A1} and~\Ref{tm-defs} in turn yields 
\bea
 && \var_0\{Z_j(i)\} \Eq  \skk \EE_0 X_k(i-1) \te\uj\uT\tSk \te\uj \non\\
  &&\qquad\quad\Le c_0 \|X(0)\|_1 \t^{i-1} \skk \m_k m_j p_{kj}(1-p_{kj})
   \Le c_{10}\|X(0)\|_1 \z\t^{i} \tm_j, 
\ena
with $c_{10} := c_0/\t$. \ep

\bigskip
We now, for future use, define the quantity 
\eq\label{zeta-star}
  \z_*\ :=\ J\max_{1\le j\le J}\|PN_Y \te\uj\|_\infty,
\en
noting also that 
\eq
   \bigl(J^{-1}\min_k \m_k\bigr)\z_* \Le \z \Eq \sjj \m^TPN_Y\te\uj \Le \|\m\|_1 \z_* 
      \Eq \z_*. \label{z-star-2}
\en
Hence, in view of~\Ref{z-reln}, we write
\eq\label{Zstar-def}
      \z_* \Def \ZZZs\sqrt{\frac mn},
\en
noting that $\ZZZ$ and~$\ZZZs$ can be thought of as having comparable magnitude.
We also introduce the notation
\eq\label{umn-def}
  u_{mn} \Def \Bl\frac mn \Br^{1/4}\Blb 1 + \Bl\frac mn \Br^{1/4} \Brb.
\en

Using
Lemma~\ref{Y-bounds}, we can now establish an analogue of 
Corollary~\ref{product-moment} for the elements of~$Y(i)$.

\begin{corollary}\label{product-moment-Y}
For all $1\le j,l \le J$, there is a constant~$c_6'$ such that
\[
  \EE_0 \{Y_j(i)Y_l(i)\} 
	  \Le c_5 \z^2\|X(0)\|_1^2 \tm_j\tm_l \t^{2(i-1)} + 
		               c_6'  u_{mn}^2 \|X(0)\|_1  i\t^{i-1}\g^{(i-1)/2},
\]
where $c_5$ is as in Corollary~\ref{product-moment}.
In particular, for $c_5'' := \ZZZ^2c_5 + c_6'J^2 \th$, 
we have
\[
  \EE_0 \|Y(i)\|_1^2 \Le c_5''  u_{mn}^2 \t^{2(i-1)} \|X(0)\|_1^2.
\]
\end{corollary}

\proof
We first note that
\bea
  \EE_0 \{Y_j(i)Y_l(i)\} &\le& |\EE_0 Y_j(i) \EE_0 Y_l(i)| + |\cov_0(Y_j(i),Y_l(i))| \non\\ 
  &\le& |\EE_0 Y_j(i) \EE_0 Y_l(i)| + \sqrt{\var_0\{Y_j(i)\}\,\var_0\{Y_l(i)\}}. \non
\ena	
Now, writing $Y_j(i) = X^T(i-1)PN_Y\te\uj + Z_j(i)$, it follows 
from Lemmas \ref{orthog-varce}\,(ii) and~\ref{Y-bounds} and from~\Ref{zeta-star} that
\beas
	\lefteqn{\sqrt{\var_0\{Y_j(i)\}}} \\
	 &&\Le \sqrt{2\|X(0)\|_1}\Blb \sqrt{c_2\t/(\t-1)}\t^{i-1}\z\tm_j
	     + \sqrt{c_3 (i-1)\g^{i-1}}J^{-1}\z_* + \sqrt{c_{10}\z\t^i/2} \Brb \\
	&&\Le \sqrt{2\|X(0)\|_1}\Blb \sqrt{c_2\t/(\t-1)}\t^{i-1}\z\tm_j
	     + u_{mn}\sqrt{c_4 i\g^{i-1}} \Brb,
\enas
where  $\sqrt{c_4} = (J^{-1}\sqrt{c_3}\ZZZs \vee \sqrt{c_{10}\ZZZ\t/2})$.  
A similar bound holds also for
$\sqrt{\var_0\{Y_l(i)\}}$.  The corollary now follows from~\Ref{A1},
and by taking $c_6' = 4 \ZZZ\sqrt{c_2c_4\t/(\t-1)} + 2c_4$.  \ep	

\subsection{$X$- and $Y$-approximation}\label{XY-approx}
Using the preparation above, we are now able to approximate $X_k(i)$
and~$Y_j(i)$,
$i\ge1$, in terms of the limiting random variable~$W$, making precise
the description at the end of Section~\ref{IB}, and bounding the error
in the approximation.  We begin by considering the $X$-components.

\begin{theorem}\label{TX-W-bnd}
There is a constant~$c_9$ such that, for $1\le k\le K$,
\eq\label{TX-appxn}
   \EE_0 \Blm X_k(i) - \t^{i}W\m_k \Brm
     \Le c_9 \|X(0)\|_1 ((i+1)\g^i)^{1/2},\qquad i\ge0.
\en
\end{theorem}

\proof
Here, we note from Lemma~\ref{orthog-varce}\,(iii), \Ref{A-2}, 
\Ref{X-mean} and~\Ref{p7star} that
\beas
   \lefteqn{\EE_0 \Blm X_k(i) - (X^T(i) \n) \m_k\Brm \phantom{HHHHHHHHHHHHH}}\\ 
     &\le&  \bigl(\var_0\{X_k(i)-(X^T(i)\n)\m_k\}\bigr)^{1/2}
         + |\EE_0\{X_k(i)-(X^T(i)\n)\m_k\}| \\
     &\le& (c_3 \|X(0)\|_1 i\g^i)^{1/2} 
         + \|X(0)\|_1\|(M_X^i - \t^i\n\m^T)e\uk\|_\infty  \\
     &\le&  c_7 \|X(0)\|_1 ((i+1)\g^i)^{1/2},
\enas
with~$c_7 := \sqrt{c_3} + c_1$,
whereas, from Lemma~\ref{mart-varce}\,(i), and the Cauchy-Schwarz inequality, 
\eq\label{var-8}
   \EE_0 \Blm (X^T(i) \n - \t^i W)\m_k\Brm
     \Le \m_k \t^i\{c_2 \{\t/(\t-1)\} \|X(0)\|_1 \t^{-i}\}^{1/2}
     \Le c_8  \|X(0)\|_1 \t^{i/2},
\en
with $c_8 := \{c_2 \t/(\t-1)\}^{1/2}$.
Hence the theorem follows, with $c_9\ :=\ c_7 + c_8$.
\ep
		 		  
\bigskip
With the help of Lemma~\ref{Y-bounds}, we can also approximate
$Y(i)$ in terms of the limiting random variable~$W$, 
complementing Theorem~\ref{TX-W-bnd}.  	

\begin{theorem}\label{TY-W-bnd}
There is a constant~$c_{14}$ such that, for each $1\le j\le J$,
\eq\label{TY-appxn}
  \EE_0 \Blm Y_j(i) - \t^{i-1}\z W\tm_j \Brm
     \Le c_{14}  u_{mn}\|X(0)\|_1 (i\g^i)^{1/2},\qquad i\ge1,
\en
where $\tm$ is given in~\Ref{tm-defs}. 
\end{theorem}

\proof
It is immediate from Lemma~\ref{Y-bounds} that
\[
   \EE_0\Blm Y_j(i) -  X^T(i-1)PN_Y \te\uj\Brm
    \Le  \{c_{10}\z\|X(0)\|_1 \|\tm\|_\infty\t^{i}\}^{1/2} 
     \Le c_{11}\z^{1/2}\|X(0)\|_1 \t^{(i-1)/2},
\]
with $c_{11} := \sqrt{c_{10}\| \tm \|_\infty}$,
and then, as in the proof of Theorem~\ref{TX-W-bnd}, using 
Lemma~\ref{mart-varce}\,(i), \Ref{p7star} and~\Ref{A-2}, we have
\beas
   \lefteqn{\EE_0\Blm X^T(i-1)PN_Y \te\uj - \z X^T(i-1)\n\tm_j \Brm}\\
    &&\Le \Bigl\{(c_3 \|X(0)\|_1 (i-1)\g^{i-1})^{1/2} 
       + c_1\|X(0)\|_1 \g^{(i-1)/2}\Bigr\}\|PN_Y \te\uj\|_\infty \\
    &&\Le c_{12} \z_*\|X(0)\|_1 (i\g^i)^{1/2},
\enas
with~$c_{12}  := J^{-1}\sqrt{1/\g}(\sqrt{c_3} + c_1) $.  
Then, once again invoking Lemma~\ref{mart-varce}\,(i), as for~\Ref{var-8},
we have
\beas
   \EE_0 \Blm \z X^T(i-1)\n\tm_j - \t^{i-1}\z W \tm_j\Brm
    &\le& c_8 \z\tm_j \|X(0)\|_1 \t^{(i-1)/2} 
    \Le c_{13}\z \|X(0)\|_1 \t^{i/2},
\enas
with $c_{13} := (c_8/\sqrt\t)\|\tm\|_\infty$. 
The theorem follows, since $\g\ge\t$, by taking
$c_{14} :=  (\{c_{11}\sqrt\ZZZ\}\vee \{ c_{12}\ZZZs +  c_{13}\ZZZ\})$. 
\ep

\section{Ghosts}\label{Ghosts}
\setcounter{equation}{0}
In the previous section, we justified simple approximations to the
joint counts $X(i)$ and~$Y(i)$ in the bipartite branching process.
We now need to show that the same approximation can be used for the
composition of the neighbourhoods in the intersection graph, albeit
with a further error. This involves showing that the effect of the
`ghosts' is not too large.
Let $G^X_k(i)$ and $G^Y_j(i)$ denote the total numbers of type $(k,1)$ 
and of type $(j,2)$ individuals of class~$0$ (ghosts), respectively, 
alive in generations~$2i$ and $2i-1$ respectively of the
bipartite process that starts with individuals $A$ and~$B$. Then
it turns out to be
enough to derive bounds for their expectations, as functions of~$i$.  

To state the result, define 
\eq\label{rho-def}
    \r^X \ :=\ \max_{1\le k\le K} {\m_k/q^X_k};\qquad 
		\r^Y \ :=\ \max_{1\le j\le J} {\tm_j/q^Y_j},
\en
where $q^X$ and~$q^Y$ are as in~\Ref{q-defs}.
Note that, for the $Z$-process under consideration, $\|X(0)\|_1 = 2$.

\begin{theorem}\label{ghost-means}
 There are constants $c^*_{15}$ and~$c^*_{16}$ such that 
\eqs
  \EE_0 G^X_k(i) &\le& c_{15}^*\t^{2i}\erra^4,  \quad 1\le k\le K;\\
  \EE_0 G^Y_j(i) &\le&  c_{16}^*\sqrt{\frac mn}\,\t^{2(i-1)}\erra^4,
        \quad 1\le j\le J,
\ens
where $\erra := n^{-1/4} + m^{-1/4} $.
\end{theorem}

\proof
The ghosts can be counted by descent from original ghosts, whose parents
were of class~$1$.  We write $H^X_k(l)$ and $H^Y_j(l)$ to denote
the numbers of original ghosts of the corresponding types born in
the $l$-th generations of the $X$ and~$Y$ processes. For $i>l\ge0$, 
we then let $\gkxyli$ and~$\gkxxli$ denote the total numbers of 
descendants of generation~$l$ original $Y$ and~$X$ ghosts, respectively, 
alive at time~$2i$ in the 
bipartite process, that are $(k,1)$ individuals; the quantities 
$\gjyyli$ and~$\gjyxli$ are defined analogously.  Thus
\eqa
   G^X_k(i) &=& H^X_k(i)  + \sum_{l=0}^{i-1} \{\gkxylpii + \gkxxli\};
	     \label{G-X-sum}\\  
   G^Y_j(i) &=& H^Y_j(i) + \sum_{l=1}^{i-1} \{ \gjyyli + \gjyxli\} 
           + \gjyxoi.  \label{G-Y-sum}
\ena  
Note that, for $i\ge l\ge0$, from~\Ref{A0} and~\Ref{z-reln},
\eq
  \EE(\gkxyli \giv H^Y(l)) 
	  \Le c_0 \ZZZ^{-1}\t\|H^Y(l)\|_1 \t^{i-l}\m_k\snm;
		    \label{B0XY}
\en
and that, for $i>l\ge0$,
\eqa
  \EE(\gkxxli \giv H^X(l)) 
	  &\le& c_0 \|H^X(l)\|_1 \t^{i-l}\m_k;
      \label{B0XX} \\
  \EE(\gjyyli \giv H^Y(l)) 
	  &\le& c_0 \|H^Y(l)\|_1 \t^{i-l}\tm_j;
     \label{B0YY} \\
  \EE(\gjyxli \giv H^X(l)) 
	  &\le& c_0 \ZZZ\|H^X(l)\|_1 \t^{i-l-1}\tm_j\smn,
     \label{B0YX} 	  
\ena
from \Ref{A0} and~\Ref{A01}.

Now an original ghost of type $(j,2)$ is created when an index from the
set $\{1,2,\ldots,m_j\}$ is re-used.  Hence
\[
    \EE\{H^Y_j(l) \giv \ff_l^Y\} 
       \Le Y_j(l)\,\Bigl\{m_j^{-1}\sum_{s=1}^l Y_j(s)\Bigr\}.
\]
Furthermore, from~\Ref{A01}, for $l > s\ge0$,
\[
    \EE\{Y_j(l) \giv \ff_s^Y\}	\Le c_0\|Y(s)\|_1 \t^{l-s}\tm_j.
\]
Combining these bounds, it follows that
\eq\label{HYj-mean}
    \EE_0 H_j^Y(l) \Le c_0 m_j^{-1}\tm_j \sum_{s=1}^l \t^{l-s}\EE_0(Y_j(s)\|Y(s)\|_1),
\en
and hence, using Corollary~\ref{product-moment-Y}, that
\eq\label{HY-mean}
    \EE_0 \|H^Y(l)\|_1 
       \Le c_0 m^{-1}\r^Y \sum_{s=1}^l \t^{l-s}\EE_0\|Y(s)\|_1^2
       \Le c_{17}m^{-1}\,\r^Y u_{mn}^2\t^{2(l-1)},
\en
where $c_{17} := 4c_5''c_0\t/(\t-1)$.
Similar calculations show that
\eq\label{HXk-mean}
    \EE_0 H_k^X(l) \Le c_0 n_k^{-1}\m_k \sum_{s=0}^l \t^{l-s}\EE_0(X_k(s)\|X(s)\|_1),
\en
and, with Corollary~\ref{product-moment}, 
\eq\label{HX-mean}
    \EE_0 \|H^X(l)\|_1 
	\Le c_{18}n^{-1}\r^X\t^{2l},
\en			 
with $c_{18} := 4c_5'c_0\t/(\t-1)$.  For future reference, we note also that,
in consequence, for any $s\ge0$,
\eqa
   \sum_{l=1}^s \t^{s-l}\EE_0 \|H^Y(l)\|_1 
      &\le& \frac{c_{17}\t}{\t-1}m^{-1}\r^Y u_{mn}^2 \t^{2(s-1)};\quad
   \non\\
   \sum_{l=0}^s \t^{s-l}\EE_0 \|H^X(l)\|_1 
      &\le& \frac{c_{18}\t}{\t-1}n^{-1}\r^X\t^{2s}. \label{H-sums}
\ena

It now remains to take expectations in \Ref{G-X-sum} and~\Ref{G-Y-sum},
using \Ref{B0XY}--\Ref{B0YX} and~\Ref{H-sums}.  For instance, for
$\EE_0 G_j^Y(i)$, we have
\eqs
  \sum_{l=1}^{i-1} \EE_0\gjyyli 
     &\le& c_0 \tm_j \sum_{l=1}^{i-1}\t^{i-l} \EE_0 \|H^Y(l)\|_1 
     \Le c_{19}m^{-1}\r^Y\tm_j  u_{mn}^2\t^{2(i-1)},
\ens
with $c_{19} := c_0 c_{17}/(\t-1)$; similar calculations yield
\eqs
  \sum_{l=0}^{i-1} \EE_0\gjyxli 
     &\le& c_{20}n^{-1}\r^X\tm_j \t^{2(i-1)}\smn,
\ens
with $c_{20} := \ZZZ c_0 c_{18}\t/(\t-1)$, and
\eqs
   \EE_0 H_j^Y(i) 
     &\le& c_{17}m^{-1}\r^Y u_{mn}^2\t^{2(i-1)},
\ens
from~\Ref{HY-mean}. Hence it follows that 
\[
   \EE_0 G^Y_j(i) \Le  \smn\,\t^{2(i-1)}\left(\frac{c_{16}\r^X}{n} 
       + \frac{\tc_{16}\r^Y}{m}\Blb 1 + \Bl \frac mn \Br^{1/4}\Brb^2 \right),
\]
with $c_{16} := c_{20}\|\tm\|_\infty$ and 
$\tc_{16} := c_{19}\|\tm\|_\infty + c_{17}$.  A similar argument yields the
bound
\[
   \EE_0 G^X_k(i) \Le  \t^{2i}\left(\frac{c_{15}\r^X}{n} 
       + \frac{\tc_{15}\r^Y}{m} \Blb 1 + \Bl \frac mn \Br^{1/4}\Brb^2 \right),
\] 
with 
\eqs
   c_{15} &:=& c_{18}\Bigl\{c_0(\t-1)^{-1}\,\|\m\|_\infty + 1\Bigr\};\quad
   \tc_{15} \ :=\ c_0 c_{17}\ZZZ^{-1}(\t-1)^{-1}\, \|\m\|_\infty.
\ens
Recalling the definition~\Ref{Zstar-def} of~$\ZZZs$, this
completes the theorem, with 
\[
     c_{15}^* \ :=\ 2(c_{15}\r^X\vee \tc_{15}\r^Y);\qquad
     c_{16}^* \ :=\ 2(c_{16}\r^X\vee \tc_{16}\r^Y).
\]
\ep

\section{The probability of a common label} \label{intersection}
 \setcounter{equation}{0}
Our next step is to establish a Poisson approximation for the probability
of a coincidence in a random labelling problem.  The underlying idea is to
look at neighbourhoods of radii $i_A$ and~$i_B$ of two initial vertices 
$A$ and~$B$; if they have no vertices in common, then the distance between
$A$ and~$B$ exceeds $i_A+i_B$. Whether two vertices in the neighbourhoods
are the same can be thought of as a labelling problem, where the
assignment of labels is almost uniform and at random.
The result that we need is the following variant
on the Poisson approximation to hypergeometric sampling.

\begin{proposition} \label{poisson1} 
For each~$l$, $1\le l\le L$, and for each~$r$, $1\le r\le R_l$, we
independently draw a subset $\{J_{lr1},\ldots,J_{lrz_{lr}}\}$ of 
size~$z_{lr}$ from a fixed set~$W_l$ of size~$w_l\ge2$, with any
subset equally likely to be drawn; define $z_l := \sum_{r=1}^{R_l}z_{lr}$.  
We then repeat the experiment
independently, with subsets $\{J'_{lr1},\ldots,J'_{lrz'_{lr}}\}$
of sizes~$z'_{lr}$, $1\le l\le L$, $1\le r\le R'_l$. For each
$w\in W_l$ and for each $l,r,s,r',s'$, define
\[
    H(w,l;r,s;r',s') \ :=\ I[J_{lrs} = J'_{lr's'} = w].
\]
Then, for fixed subsets $W^*_l \subset W_l$ of sizes~$w^*_l$, 
$1\le l\le L$, define
\[
   S\ :=\ S(z,z';w,w^*)\ :=\ 
    \sll\swls\sum_{r=1}^{R_l}\sum_{r'=1}^{R'_l}\sum_{s=1}^{z_{lr}}  
     \sum_{s'=1}^{z'_{lr'}} H(w,l;r,s;r',s')
\]
to be the number of pairs of elements, one from each sample, that 
consist of two copies of the same element, not
belonging to any of the~$W^*_l$.  Then
\beas
   \left|\PP[S(z,z';w,w^*)=0] - \exp({-\lambda(z,z',w,L)}) \right| 
	      &\leq& B_1(z,z',w,L) + B_1^*(z,z',w,w^*,L), 
\enas 
where
\bea
    \l(z,z',w,L) &:=& \sll \frac{z_l z'_l}{w_l}; \qquad 
    B_1(z,z',w,L) \ :=\ 2\sum_{l=1}^L \frac{z_l+z_l'}{w_l} , \label{B1-def}
\ena
and
\eqa    
    B_1^*(z,z',w,w^*,L) &:=& \sll\frac{z_l z'_l w^*_l}{w_l^2}. \label{B1star-def}
\ena
\end{proposition} 

\msk 
\proof
The indicators $H(w,l;r,s;r',s')$ are negatively related (for the definition see~\cite{bhj} Definition~2.1.1), as can be seen
by constructing an explicit coupling very much as in~\cite{bhj}, p.112. 
Setting
\[
   \tH(l,r,r')\ :=\ 
     \swls \sum_{s=1}^{z_{lr}} \sum_{s'=1}^{z'_{lr'}} H(w,l;r,s;r',s'),
\]
the random variables $\tH(l,r,r')$ are pairwise independent, and satisfy
\beas
   &&\EE\tH(l,r,r') 
     \Eq \left(1 - \frac{w_l^*}{w_l}\right)\frac{z_{lr}z'_{lr'}}{w_l};\\
   &&0 \Le \EE\tH(l,r,r') - \var\tH(l,r,r') 
       \Le \EE\tH(l,r,r') \left(\frac{z_{lr} + z'_{lr'} - 1}{w_l-1}\right);\\
   && S \Eq \sll\sum_{r=1}^{R_l}\sum_{r'=1}^{R'_l} \tH(l,r,r').
\enas  
Hence, since $w_l\ge2$ for all~$l$, it follows that
\beas
   0 &\le& \EE S - {\var S}
     \Le 2\max_{l,r,r'}\left(\frac{z_{lr} + z'_{lr'}}{w_l}\right) \EE S,
\enas
so that
\beas 1 - \frac{\var S}{\EE S}
   &\le& 2\sll \sum_{r=1}^{R_l}\frac{z_{lr}}{w_l} 
           + 2\sll\sum_{r'=1}^{R'_l} \frac{z'_{lr'}}{w_l}
             \Le B_1(z,z',w,L).
\enas
From this, and since
\[
   \sll\sum_{r=1}^{R_l}\sum_{r'=1}^{R'_l}\EE\tH(l,r,r') 
     \Eq \lambda(z,z',w,L) - B_1^*(z,z',w,w^*,L),
\]
it follows using~\cite{bhj} Theorem~2.C.2 that
\eq\label{appx-1}
   \bigl|\PP[S(z,z';w,w^*)=0] 
         - \exp\{{-[\lambda(z,z',w,L) - B_1^*(z,z',w,w^*,L)]}\}\bigr|
    \Le B_1(z,z',w,L).
\en
The final estimate follows because 
\eq\label{etaylor}
|e^{-\l} - e^{-\l'}| \le \min\{1,|\l-\l'|\}
\en  when
$\l,\l' \ge 0$.
\ep

\medskip
\nin Note that, if, for some $\t > 1$, none of $z_l$ and $z'_l$
exceeds $C_1\t^{i/2}$ and all of the $w_l$ exceed $C_2\t^i$, then
\beas 
   B_1(z,z',w,L) &\leq & \frac{4LC_1}{C_2\,\t^{i/2}}, 
\enas 
geometrically small with~$i$.
     
Now suppose that we do not know the true values $z_l$ and $z'_l$,
but only approximations  $\tz_l$ and $\tz'_l$ to them.  
Then we can instead use these to approximate $\PP[S(z,z',w,w^*)=0]$,  with 
some possible extra error.

\begin{proposition} \label{poisson2}
Suppose that 
\[
   |\tz_l - z_l| = \e_l
	    \quad \mbox{and}\quad |\tz'_l - z'_l| = \e'_l,\qquad 1 \le l \le L.
\]			
Then
\beas
   \lefteqn{\left|\PP[S(z,z',w,w^*)=0] - \exp({-\lambda(\tz, \tz', w,L)}) \right|}\\ 
   &&\Le B_1(z,z',w,L) + B_1^*(z,z',w,w^*,L) + B_2(z,z',\e,\e',w,L), 
\enas 
where
\beas 
    B_2(z,z',\e,\e',w,L) &=& 
		   \lb(z,\e',w,L) + \lb(\e,z',w,L) + \lb(\e,\e',w,L),
\enas 
where $\lb := \min(\l,1)$. 
\end{proposition} 

\proof
Immediate from \Ref{etaylor}.
\ep

In Section~\ref{distances}, we take for $\tz$ and~$\tz'$ convenient approximations 
to numbers of individuals alive in generations~$r\ge1$ in the bipartite 
branching process that starts from the two individuals~$A$ of type 
$(k_1,1)$ and~$B$ of type~$(k_2,1)$.  For $r=2l$ even, we take~$\tz$
to approximate the~$X^A(l)$ descendants of~$A$, and~$\tz'$ to approximate
the~$X^B(l)$ descendants of~$B$, and $w_l = n_l$, $1\le l\le K$. 
For $r=2l-1$ odd, we take~$\tz$ to approximate the~$Y^A(l)$ descendants 
of~$A$, and~$\tz'$ to approximate the~$Y^B(l)$ descendants of~$B$, now
with $w_l = m_l$, $1\le l\le J$. We then show that these approximations 
are sufficiently close to the corresponding numbers $z$ and~$z'$
of class~$1$ and class~$0'$ descendants of individuals $A$ and~$B$
in generation~$r$, 
so that $\PP[S(z,z',w) = 0]$ is correspondingly close to 
$\exp\{-\l(\tz,\tz',w)\}$.

\section{Approximating inter-point distances}\label{distances}
\setcounter{equation}{0}
We now return to the problem of real interest, the distribution
of the graph distance $D := D_{k_1,k_2}$ between two vertices 
$(k_1,1)$ and $(k_2,1)$ in the intersection graph, taken to
be infinite if the vertices are in different components. 

\subsection{Conditioning on the branching process}\label{bp-condition}
We begin by approximating the conditional probability 
$\PP[D > d\giv Z]$ that $A$ and~$B$ are more than distance~$d$ apart, 
given the trajectory of
the bipartite process~$Z$ starting from $A$ and~$B$.  The
conditional probability is then a function only of the way 
in which the labels were assigned to the individuals in the
process~$Z$. The labelling determines the classes of the
individuals, and the event $\{D>d\}$ occurs exactly when there are no overlaps
between the labels of the class~$1$ and class~$0'$ individuals that are 
descended from $A$ and those of the descendants of~$B$, at any generation~$l$, 
$1\le l\le d$, of~$Z$.  Let $\ggg_s$ denote the information in the 
labels up to generation~$s$. 

\begin{proposition}\label{conditional-prob}
For any $1\le l\le d$,
\eq
   \bigl|\PP[D > l \giv Z,\ggg_{l-1} \cap \{D \ge l\}] 
       - e^{-\l'(l,Z)}\bigr| \Le \f(l,Z),
\en
where~$\f$ is given in \Ref{total-error-even} and~\Ref{total-error-odd},
and $\l'(l,Z)$ in \Ref{ldash-even} and~\Ref{ldash-odd}.
It then follows in particular that 
\eq\label{prod-form}
   \left|\PP[D > d \giv Z] - \prod_{l=1}^{d} e^{-\l'(l,Z)}\right| \Le 
    \sum_{l=1}^{d} \f(l,Z).
\en
\end{proposition}

\proof
 Suppose, first, that~$l$ is even.
By Proposition~\ref{poisson1}, the probability
\[
    \PP[D > l \giv Z,\ggg_{l-1} \cap \{D \ge l\}]
\]
is close to $\exp\{-\l(z(l),z'(l),\tnn,K)\}$, with $\tnn_k=n_k$,
$z(l)$ the numbers of children of the different types of class~$1$ 
descendants of~$A$ in generation~$l-1$, and $z'(l)$ is the same 
for descendants of~$B$. When applying Proposition~\ref{poisson1},
$R_l$ represents the number of class~$1$ descendants of~$A$ in generation~$l-1$,
and~$z_{lr}$ the number of offspring of the $r$-th of these;
these offspring make up the class~$1$ and class~$0'$ descendants 
of~$A$ in generation~$l$. 
The error in the approximation is then no larger than 
$$
   B_1(z(l),z'(l),\tnn,K) + B_1^*(z(l),z'(l),\tnn,\tnn^*(l-1),K),
$$   
a quantity that we shall need to bound later, 
where $\tnn^*_k(l-1)$ is the number of labels for $(k,1)$ individuals 
already used up to  generation~$l-1$ of the $Z$-process.  

The quantities $z(l)$ and~$z'(l)$ appearing in~$\l$ are not directly 
accessible, and are not functions of~$Z$ alone.
However, we can exploit Proposition~\ref{poisson2}, provided that
we can find suitable approximations to them.  The first is to
replace $z(l)$ by~$X^A(l/2)$, the numbers of descendants of~$A$
of the different types in generation~$l/2$ of the $X$-process, and 
$z'(l)$ by~$X^B(l/2)$,
noting that
\eq\label{G-bnds}
   0 \Le X_k^A(l/2) - z_k(l) \Le G_k^X(l/2); \quad
   0 \Le X_k^B(l/2) - z'_k(l) \Le G_k^X(l/2),
\en
and that the number of ghosts~$G_k^X(l/2)$, investigated in 
Section~\ref{Ghosts}, is calculated for the whole bivariate
process.  Then the quantities $X^A(l/2)$ and~$X^B(l/2)$ can in
turn be more simply approximated, using Theorem~\ref{TX-W-bnd},
by $\t^{l/2}W^A\m$ and $\t^{l/2}W^B\m$,
where $W^A$ is the limit of the martingale $\t^{-i}\n^T X^A(i)$,
and~$W^B$ the limit of $\t^{-i}\n^T X^B(i)$: note that these two
random variables are independent, by the branching property.
From Proposition~\Ref{poisson2}, replacing $z(l)$ by $\t^{l/2}W^A\m$ and $z'(l)$ by 
$\t^{l/2}W^B\m$, we incur a further error of at most 
$B_2(X^A(l/2),X^B(l/2),\e^A(l),\e^B(l),\tnn,K)$, where
\eq\label{epsilons}
  \e^A_k(l) \Eq G_k^X(l/2) + E(l/2,X^A); \qquad
  \e^B_k(l) \Eq G_k^X(l/2) + E(l/2,X^B),
\en
with $E(l/2,X^A) = X^A(l/2) - \tau^{l/2} W^A \mu$ and  $E(l/2,X^A)$ defined analogously. By Theorem~\ref{TX-W-bnd}, 
\eq\label{E-bnd}
   \EE\{ E(i,X) \} \Le c_9((i+1)\g^{i})^{1/2},
\en
for $X=X^A$ and for $X=X^B$.   

We also clearly have $\tnn^*(l-1) \le T^X((l-2)/2)$,
componentwise, where
\[
    T^X(s)\ :=\ \sum_{r=0}^s X(r) = \sum_{r=0}^s  (X^A(r) + X^B(r)),
\]
and, as observed above, $z(l) \le X^A(l/2)$, $z'(l) \le X^B(l/2)$.  As a result,
the approximation error at this step is no larger than
\eqa
  \f(l,Z) &:=& B_1(X^A(l/2),X^B(l/2),\tnn,K) 
       + B_1^*(X^A(l/2),X^B(l/2),\tnn,T^X((l-2)/2),K)  \non\\
   &&\quad\mbox{} + B_2(X^A(l/2),X^B(l/2),\e^A(l),\e^B(l),\tnn,K), 
       \label{total-error-even}
\ena
with $\e^A(l),\e^B(l)$ as in~\Ref{epsilons};
thus we have, for $l$ even,
\eq\label{one-step-X}
  |\PP[D > l \giv Z,\ggg_{l-1}\cap \{D \ge l\}] - \exp\{-\l'(l,Z)\}|
    \Le \f(l,Z),
\en
where, for~$l$ even, 
\eq\label{ldash-even}
  \l'(l,Z) \ :=\ \l(\t^{l/2}W^A\m,\t^{l/2}W^B\m,\tnn,K).
\en
A similar argument for~$l$ odd yields    
\eq\label{one-step-Y}
|\PP[D > l \giv Z,\ggg_{l-1}\cap \{D \ge l\}] - \exp\{-\l'(l,Z)\}|
    \Le \f(l,Z),
\en
where, for $l$ odd,
\eqa
   \l'(l,Z) &:=& \l(\z\t^{(l-1)/2}W^A\tm,\z\t^{(l-1)/2}W^B\tm,\tmm,J),
     \label{ldash-odd} \\
   \f(l,Z) &:=& B_1(Y^A((l+1)/2),Y^B((l+1)/2),\tmm,J) \non\\
      &&\quad\mbox{} + B_1^*(Y^A((l+1)/2),Y^B((l+1)/2),\tmm,T^Y((l-1)/2),J) \non\\
       &&\qquad\mbox{} + B_2(Y^A((l+1)/2),Y^B((l+1)/2),\h^A(l),\h^B(l),\tmm,J).
          \label{total-error-odd}
\ena 
Here,
\eq\label{etas}
  \h^A_k(l) \Eq G_k^Y((l+1)/2) + E'((l+1)/2,Y^A); \qquad
  \h^B_k(l) \Eq G_k^Y((l+1)/2) + E'((l+1)/2,Y^B),
\en
and
$$
   \EE\{ E'(i,Y) \} \Le c_{14}u_{mn}(i\g^{i})^{1/2},
$$
for $Y=Y^A$ and for $Y=Y^B$, by Theorem~\ref{TY-W-bnd}.  
This proves the first statement of the proposition. 

The second part is easier.  We first note that
\[
   \PP[D > l \giv Z] \Eq \EE\{\EE(I[D > l-1] 
     \PP[D > l \giv Z,\ggg_{l-1}\cap \{D \ge l\}]\giv \ggg_{l-1},Z)\giv Z\},
\]
and deduce from the first part that
\[
   \bigl|\PP[D > l \giv Z] - e^{-\l'(l,Z)}\PP[D > l-1 \giv Z]\bigr|
     \Le \f(l,Z),
\]
from which the last part follows.
\ep

\bsk
Thus, combining \Ref{one-step-X} and~\Ref{one-step-Y} with Proposition~\ref{conditional-prob}, we
find that
\eqa\label{all-steps}
   | \PP[D > d \giv Z] - \exp\{-W^AW^B L(d)\} | &\le& \sum_{l=1}^{d} \f(l,Z) ,
\ena   
where, using~\Ref{z-reln},
\eqs
   L(2i) &:=& \frac{\t^{2i}-1}{\t^2-1} \Blb \t^2\skk  \frac{\m_k^2}{n_k} 
             + \z^2\sjj \frac{\tm_j^2}{m_j} \Brb\\
        &=& \frac{\t^{2i}-1}{\t^2-1}\, \t( \t + 1)   \skk  \frac{\m_k^2}{n_k} 
            \Eq \Bl\frac{\t}{\t-1}\Br n^{-1}(\t^{2i}-1) \skk  \frac{\m_k^2}{q_k^X};\\
   L(2i+1) &:=& \frac{1}{\t^2-1} \Blb (\tau^{2i+2} - \tau^2) \skk  \frac{\m_k^2}{n_k} 
             + (\tau^{2i+2} - 1) \z^2\sjj \frac{\tm_j^2}{m_j} \Brb \\
   &=& \frac{1}{\t^2-1} \Blb ( \tau^{2i+2} - \tau^2 +  \tau^{2i+3} - \t) \skk 
       \frac{\m_k^2}{n_k}   \Brb \\
   &=& \Bl\frac{\t}{\t-1}\Br n^{-1}(\t^{2i+1}-1) \skk  \frac{\m_k^2}{q_k^X},
\ens  
and hence, with~\Ref{kappa-def}, 
\eq
    L(d) \Eq \kappa n^{-1}(\t^{d}-1),
\en
for~$d$ both even and odd.

\subsection{The unconditional distribution}\label{D-distn}
The unconditional probabilities for~$D$ are now given by taking expectations
in conjunction with~\Ref{all-steps}, so that it just remains to evaluate
the terms $\EE_{k_1,k_2}\f(l,Z)$.  Note that, in the approximation,
randomness comes in only through the independent random variables
$W^A$ and~$W^B$, the first with a distribution which depends only on the
value of~$k_1$, and the second on~$k_2$.

To assist in judging the impact of the various factors in our
bounds, it is convenient to define 
\eq\label{i0-def}
   i_0 \ :=\ \left\lfloor \frac{\log n}{\log\t} \right\rfloor,
\en
so that $\t^{i_0} \le n < \t^{i_0+1}$, and to set
\beas
   \t^{-1} \ <\ \phh(n) \Def n^{-1}\t^{i_0} \Le 1.
\enas  
Then, with $\kappa$ given in~\Ref{kappa-def} and under standard asymptotics, 
\eq\label{L-approx}
   |L(d) - \t^{d-i_0}\kkk\phh(n)| \Le \frac{\r^X}{n} \Bl\frac{\t}{\t-1}\Br\ \to\ 0,
\en
and  $\kkk$ remains bounded away from $0$ and~$\infty$. 

\begin{theorem}\label{ip-distance}
For $d = i_0 + u$, with $u\in\integ$ and $|u| < i_0/2$, we have
\beas
  \lefteqn{|\Pkk[D - i_0 > u] - \Ekk\exp\{-W^AW^B\kkk\t^u\phh(n)\}|}\\
    &&\qquad\qquad\Le c_{25}\bigl\{(\t^{3u/2}+1)(n^{1/4}\erra^2\wedge1) +
         (\t^{u}+1)n^{1/4}\erra\tth_{i_0}\bigr\},  
\enas
for a suitable quantity~$c_{25}$, where 
\eq\label{theta-tilde-def}
  \tth_i \Def (i+1)^{1/2}(\g/\t^2)^{i/4}.
\en
\end{theorem}

\proof
The approximating expression is immediate, from \Ref{all-steps} and~\Ref{L-approx}, 
incurring an error of at most
\[
   \n_{k_1}\n_{k_2} \frac{\r^X}{n} \Bl\frac{\t}{\t-1}\Br.
\]
For the rest, we just need to investigate $\Ekk\f(l,Z)$ for 
$1 \le l \le i_0 + u$.

To start with, for $l = 2r$, we have
\eq\label{bound-1}
   \Ekk B_1(X^A(r),X^B(r),\tnn,K) 
           \Le 4 n^{-1} K c_0 \t^r \r^X\ \le\ n^{-1}c_{22}\t^r,
\en
with $c_{22} := 4Kc_0\r^X$, from \Ref{B1-def} and~\Ref{A1}.  Then
\eqa
   \lefteqn{\Ekk B_1^*(X^A(r),X^B(r),\tnn,T^X(r-1),K)}\non \\ 
     &&\Eq \skk \sum_{s=1}^{r-1}(nq^X_k)^{-2}
         \Ekk \{X^A_k(r)X^B_k(r)[X^A_k(s) + X^B_k(s)]\}\non  \\
   &&\Le c_{23}n^{-2}\t^{3r},  \label{bound-2}
\ena
where $c_{23} = 2K(c_0\r^X)^2(c_5\|\m\|_\infty + c_6K\th)/(\t-1)$,
from \Ref{A0} and Corollary~\ref{product-moment}.      
For $\Ekk B_2(X^A(r),X^B(r),\e^A(2r),\e^B(2r),\tnn,K)$, with
$\e^A$ and~$\e^B$ defined as in~\Ref{epsilons}, we need to be
a little more careful, because of the
product $G^X_k(r)(X_k^A(r) + X_k^B(r))$.  However, from 
Theorem~\ref{ghost-means}, it follows by Markov's inequality that,
for any $\fff = \fff(m,n,r)$,
\eqa
   \Pkk\left[\max_{1\le k\le K} 
            G_k^X(r) > \fff \right] 
   \Le \Pkk\left[\skk 
            G_k^X(r) > \fff \right]
   \Le \fff^{-1} \t^{2r}
         Kc_{15}^*\erra^4,\phantom{H} \label{Markov}
\ena
and because $B_2$ can never exceed the value~$3$, it follows that
\eqa
   \lefteqn{\Ekk B_2(X^A(r),X^B(r),\e^A(2r),\e^B(2r),\tnn,K) } \non\\
     &&\Le \Ekk B_2(X^A(r),X^B(r),\teps^A(2r),\teps^B(2r),\tnn,K) +
        6  \Pkk\left[\max_{1\le k\le K}  G_k^X(r) > \fff \right] \non\\
     &&\Le \Ekk B_2(X^A(r),X^B(r),\teps^A(2r),\teps^B(2r),\tnn,K
        + 6c_{24}  \fff^{-1}\t^{2r}\erra^4, \label{bound-3}
\ena
with 
$$
    \teps^A_k(2r) \Eq \fff + E(r,X^A)\quad\mbox{and} \quad 
          \teps^B_k(2r) \Eq \fff + E(r,X^B),
$$
and with $c_{24} := Kc_{15}^*$.  
But now, from \Ref{E-bnd} and~\Ref{A1}, it follows that
\eqa
   \lefteqn{\Ekk B_2(X^A(r),X^B(r),\teps^A(2r),\teps^B(2r),\tnn,K)}
      \label{B2-final}\\
     && \Le 2K\r^X n^{-1} c_0\t^r\{\fff + c_9((r+1)\g^r)^{1/2}\}
         + 2n^{-1}\{\fff^2 + c_9^2 (r+1)\g^r \}\skk \frac1{q_k^X}. \non
\ena
Choosing $\fff^2(m,n,r) := \t^r n \erra^4$, and then
adding the contributions from \Ref{bound-1} --~\Ref{B2-final} for 
$1\le r \le \lfloor (i_0+u)/2 \rfloor$ gives, after some computation, a bound
of the form
\eq\label{summary-1}
   c'_{25}\bigl(n^{1/4}(\t^{3u/4}+1)\erra^2 + (\t^{u}+1)\tth_{i_0}
       + (\t^{3u/2}+1)n^{-1/2}\bigr).
\en

Bounds analogous to \Ref{bound-1} and~\Ref{bound-2} hold also for $l = 2r-1$, 
with $Y,J,m$ replacing $X,K,n$ throughout the argument and estimates, and 
with $c_{22}$ and~$c_{23}$ replaced by $c'_{22}\smn$ and $c'_{23}\smn u_{mn}^2$,
where
$$
   c'_{22} \Eq 4c_0\ZZZ\t^{-1}J\r^Y \quad \mbox{and}\quad
    c'_{23} \Eq 2JZ(c_0\r^Y/\t)^2 (c_5 \ZZZ^2\t^{-2}\|\tm\|_\infty
                 + c_6'J\th)/(\t-1).
$$
The bound corresponding to~\Ref{Markov} is
\eqa
  \Pkk\left[\max_{1\le j\le J}   G_j^Y(r) > \fff \right]
   \Le \fff^{-1} \smn\,\t^{2(r-1)}
         Jc_{16}^*\erra^4\,, \label{Markov-Y}
\ena
and we also have
\eqa\label{B2-final-Y}
   \lefteqn{\Ekk B_2(Y^A(r),Y^B(r),\teta^A(2r-1),\teta^B(2r-1),\tmm,J)}\\
      &\Le & 2J\r^Y m^{-1}\t^{-1}c_0\t^r\ZZZ\smn\, \{\fff + c_{14}u_{mn}(r\g^r)^{1/2}\}
         + 2m^{-1}\{\fff^2 + (c_{14}u_{mn})^2 r\g^r \}\sjj\frac1{q_j^Y},\non
\ena
with
$$
    \teta^A_k(2r-1) \Eq \fff + E'(r,Y^A)\quad\mbox{and} \quad 
          \teta^B_k(2r-1) \Eq \fff + E'(r,Y^B).
$$
Here, we take $\fff^2 := m\t^r\erra^4$ in \Ref{Markov} and~\Ref{B2-final-Y}, and 
then, adding the errors over $1 \le r \le  \lfloor i_0+1+u \rfloor$, and
after much calculation, a bound of the form
\eq\label{summary-2}
   c''_{25}\Blb n^{1/4}\erra^2(\t^{3u/4}+1) + n^{1/4}\erra(\t^{u}+1)\tth_{i_0}
       + (\t^{3u/2}+1)\snm\erra^2\Brb
\en
is obtained.  

To deduce the bound given in the theorem, it now suffices to observe that
\[
      \bigl( n^{1/4}\erra^2 \bigr)^2 
       \ \ge\ n^{-1/2} + \snm\,\erra^2,
\]
so that the final terms in \Ref{summary-1} and~\Ref{summary-2} can be
absorbed into the first term, if the larger of the $\t$-exponents
is used. 
\ep

The defective real valued random variable~$U$, whose distribution function
\eq\label{U-dist}
  \Pkk[U \le u]  \Eq 1 - \Ekk\exp\{-W^AW^B\kkk\t^u\phh(n)\}
\en
approximates that of $D-i_0$ for integer arguments, 
can be expressed as a (defective) 
translation mixture of scaled negative standard Gumbel 
random variables. If $W^AW^B$ has distribution function~$F_{k_1,k_2}$ 
on~$\re_+$, and if
\eq\label{U-dash}
   \Pkk[U' \le u] \Def \int_{(0,\infty)} 
    \PP[-(\log\t)^{-1}(\G + \log x + \log\kkk) \le u]\, 
      dF_{k_1,k_2}(x),
\en
where~$\G$ denotes a standard Gumbel random variable, then 
\[
   \Pkk[U \le u] \Eq \Pkk[U' \le u + \log\phh(n)/\log\t].
\]
Alternatively, we can write 
\eq\label{U-infty}
   \pr[U=\infty] \Eq \pr[U'=\infty] \Eq F_{k_1,k_2}(0) \Eq 1 - \PP_{k_1}[W>0]\PP_{k_2}[W>0],
\en
and express the distribution $\law(U'\giv U'<\infty)$ as that of a random
variable~$\tU$, realized as
\eq\label{tU-def}
     \tU \Eq -\frac1{\log\t}\{\G + \log\tW_A + \log\tW_B + \log\kkk\},
\en
where $\G$, $\tW_A$ and~$\tW_B$ are independent,
\[
      \pr[\tW_A \le w] \Eq \pr_{k_1}[W \le w \giv W > 0]
   \quad\mbox{and}\quad \pr[\tW_B \le w] \Eq \pr_{k_2}[W \le w \giv W > 0].
\]
Note that $F_{k_1,k_2}(0) = \pr[U'=\infty]$ indeed approximates the 
probability that $(k_1,1)$ and $(k_2,1)$ are in different 
components of the graph, and are hence at infinite distance from 
one another, as can be seen in the following result.

\begin{theorem}\label{W=0}
There are constants $\t_1 > 1$ and $c_{26} < \infty$ such that
$$
   \max\{\Pkk[D < \infty \giv W^A = 0], \Pkk[D < \infty \giv W^B = 0]\} 
      \le c_{26} \t_1^{-i_0}.
$$   
\end{theorem}

\proof
We make the calculation for~$A$; for~$B$ the argument is the same.
From the general theory of multi-type branching processes, 
see for example \cite{daly} or \cite{JagersLageras}, conditional on the
event $\{W^A=0\}$, $X^A$ is a subcritical  branching process,
and there exist $\t_1 > 1$ and $C < \infty$ such that 
$$
    \EE_k[\|X^A(i)\|_1 \giv W^A = 0] \Le C\t_1^{-i} \quad\mbox{for all}\ k,
$$
and thus 
\[
   \EE_k[\|Y^A(i)\|_1 \giv W^A = 0] \Le C\t_1^{-i+1}c_0\z \quad\mbox{for all}\ k.
\]
Hence, immediately,     
$$
    \Pkk[2i_0 < D < \infty \giv W^A = 0] \Le C\t_1^{-i_0}.
$$
Then, for $1 \le i \le i_0$, from~\Ref{A1},
\eqs
   \Pkk[D = 2i \giv W^A = 0] &\le& 
       \skk n_k^{-1}\Ekk X^A_k(i)\Ekk X^B_k(i) 
             \Le n^{-1}KC\t_1^{-i}c_0\t^i\r^X; \\   
   \Pkk[D = 2i-1 \giv W^A = 0] &\le& m^{-1}JC\z^2\t_1^{-i+1}c_0\t^{i-1}\r^Y
     \Le n^{-1}JC\ZZZ^2 \t_1^{-i+1}c_0\t^{i-1}\r^Y,
\ens   
and the theorem follows by adding over $1 \le i\le i_0$.
\ep

In view of the considerations above, our approximation can be
summarized as follows.

\begin{corollary}\label{Main}
    For $d = i_0 + u$, with $u\in\integ$ and $|u| < i_0/2$, we have
\beas
  &&|\Pkk[D \le u + i_0] - \Pkk[U' \le u + \log\phh(n)/\log\t]| \Le
     \d(\t^u,m,n); \\
   &&|\Pkk[D = \infty] - \Pkk[U' = \infty]| \\
   &&\qquad  \Le \d(n^\a,m,n) + n^{-1}
      + \Pkk[0 < W^AW^B \le \t^2\kkk^{-1} n^{-\a}\log n] + 2c_{26}\t_1^{-i_0},
\enas
for any~$0<\a< (i_0-2)/2i_0\approx 1/2$, where
\[
   \d(y,m,n) \Def c_{25}\bigl\{(y^{3/2}+1)(n^{1/4}\erra^2\wedge1) +
         (y+1)n^{1/4}\erra\tth_{i_0}\bigr\},
\]
 $i_0$ is as in~\Ref{i0-def}, $\tth_i$ is as in~\Ref{theta-tilde-def},
$\t_1$ is as for Theorem~\ref{W=0} and~$U'$ has distribution given either 
by~\Ref{U-dash} or by \Ref{U-infty} and~\Ref{tU-def}.
\end{corollary}

\proof
The first inequality is from Theorem~\ref{ip-distance}.  For the second, we
have
\eqs
    \lefteqn{\Pkk[D < \infty]}
   \\ &&\Le \Pkk[D < \infty \giv W^A = 0] + \Pkk[D < \infty \giv W^B = 0]
        + \Pkk[W^AW^B>0],
\ens
giving
\[
   \Pkk[D = \infty] \ \ge\ 1 - \Pkk[W^AW^B>0] - 2c_{26}\t_1^{-i_0}
     \Eq \Pkk[U'=\infty] - 2c_{26}\t_1^{-i_0}.
\]
On the other hand, taking $u = \lfloor\a i_0\rfloor$ in the first part, we
have 
\[
   \Pkk[D = \infty] \Le 1 - \Pkk[U \le u] + \d(n^\a,m,n),
\]
and, from~\Ref{U-dist}, for any $C>0$,
\eqs
     \Pkk[U \le u] &\ge& 1 - \Pkk[W^AW^B \le Cn^{-\a}\log n] -
    \exp\{-Cn^{-\a}\log n \kkk (n^\a/\t) \phh(n)\}\\
    &\ge& \Pkk[W^AW^B>0] - \Pkk[0 < W^AW^B \le C n^{-\a}\log n]\\
       &&\mbox{}\hskip3in - \exp\{-C\log n \kkk /\t^2\}.
\ens
Hence, taking $C=\t^2/\kkk$, 
\eqs
   \lefteqn{\Pkk[D = \infty]}\\
   && \Le \Pkk[U'=\infty] + \Pkk[0 < W^AW^B \le \t^2\kkk^{-1} n^{-\a}\log n]
      + n^{-1} + \d(n^\a,m,n),
\ens
and the corollary is proved.
\ep

\bigskip
\nin{\bf Remark.}\  The corresponding result for the
unipartite Erd\H os--R\'enyi graph may also be of interest, 
although, as discussed at the end of Section~\ref{intro}, it is not
directly useful for our purposes.  For such a graph, the vertices are
divided into~$K$ types, with $n_k$ of type~$k$, $1\le k\le K$, and
with $n := \skk n_k$.
Edges are then independently assigned, with probabilities $p_{k,k'}$ depending
on the vertex types $k$ and~$k'$: the matrix~$P$ is thus symmetric.
The mean matrix~$M$ for the associated branching process is given by~$PN$, where
$N := {\rm diag}\{n_1,\ldots,n_K\}$, and we assume that it is irreducible
and aperiodic, and that its largest
eigenvalue $\tilt > 1$.  With these assumptions, and writing~$\m^T$ for 
the left eigenvector of~$M$ with eigenvalue~$\tilt$, only small changes need 
to be made to the sketched argument concluding Section~\ref{IB}.  
Considering coincidences in the indices
in order of increasing generation number, and with the offspring
of~$A$ considered before those of~$B$, links in the Erd\H os--Re\'nyi
graph arise exactly when there are coincidences between indices
of the~$X^A_k(i)$ and those of the~$X^B_k(i-1)$, or between indices
of the~$X^B_k(i)$ and those of the~$X^A_k(i)$.  This leads to an
approximate mean number of coincidences, up to and including the time 
when the first~$i$ generations of descendants of~$A$ and the first 
$i-1$ of~$B$ have been considered, of $\tkkk n^{-1}\tilt^{2i-1} W^A W^B$, where 
\[
     \tkkk \Def \frac\tilt{\tilt-1} \skk \frac{\m_k^2}{q_k},
\]
and $q_k := n_k/n$.
For the time until the first~$i$ generations of both have been
considered, the corresponding approximation is $\tkkk n^{-1}\tilt^{2i} W^A W^B$.
This gives the probability that the distance between $A$ and~$B$
exceeds~$d$ as being approximately
\eq\label{ER-result}
     \Ekk\Blb e^{- \tkkk n^{-1}\tilt^{d} W^A W^B} \Brb,
\en
very much the same as the formula in Theorem~\ref{ip-distance}.
Note once again that the assumption of irreducibility prevents this
line of argument being directly applicable to the bipartite model.

\subsection{Asymptotic behaviour}\label{asymptotics-2}
Recalling the standard asymptotics of Section~\ref{offspring},
we now distinguish the possibilities for the bipartite branching process
starting with a single vertex of type $(k,1)$, according to the behaviour of the
ratio $m/n$, as $\nti$. This, in turn, enables one to deduce the asymptotic
form of the approximating random variable~$U'$.
  
First, note that the~$Y_j\umn(1) \sim \Bi(m\tq_j^Y,(mn)^{-1/2}\Pi_{kj})$, 
$1\le j\le J$, are independent. 
If $m/n \to r$ with $0 < r < \infty$, then the Poisson approximation to
the binomial distribution thus shows that the distribution of~$Y_j\umn(1)$ 
differs in total variation from $\Po(\sqrt r \Pi_{kj}\tq_j^Y(m,n))$ by at
most $(mn)^{-1/2}\Pi_{kj} \sim n^{-1}\Pi_{kj}/\sqrt r$, and, conditional
on~$Y\umn(1)$, the distributions of the~$X_l\umn(1)$ are independent, and
close to the same order to $\Po(\sjj Y_j\umn(1)\Pi_{lj}\tq_l^X(m,n)r^{-1/2})$.
Hence, as $m$ and~$n$ tend to infinity in this way, the bipartite
branching process converges to the one with exactly Poisson
offspring distributions and with $Q_X\umn$ and $Q_Y\umn$ replaced
by $Q_X$ and~$Q_Y$.  Hence the distribution 
$\law(W\umn \giv X(0) = e\uk)$ converges to $\law(W \giv X(0) = e\uk)$,
for each~$k$, where~$W$ is the limiting random variable associated
with the limiting Poisson--based branching process.  It thus follows 
that the distribution of the random variable~$U'{\vphantom{X}}\umn$ also converges
to that of the corresponding~$U'$. However, the distribution
of~$U\umn$ does not converge in general, because the value of
$\log\phh(n)/\log\t$ oscillates between $-1$ and~$0$ as~$n$ varies.

If $m/n \to \infty$, the Poisson approximation 
$\Po(\sqrt {m/n} \Pi_{kj}\tq_j^Y(m,n))$ to the distribution 
of~$Y_j\umn(1)$ still has error of at most $(mn)^{-1/2}\Pi_{kj}$.
However, a simple calculation shows that, for $B_j\umn$ a Bernoulli
random vector with $\PP[B_j\umn = e\ul] = \tq_l^X(m,n)\Pi_{lj}\sqrt{n/m}$,
\[
    d_{TV}\bigl(\law(X\umn(1) \giv Y\umn(1) = \te\uj),\law(B_j\umn)\bigr) 
      \Le \Blb \skk \tq_k^X(m,n)\Pi_{kj}\sqrt{n/m} \Brb^2.
\]
It now follows from the Poisson thinning theorem (see for example Chapter 8 Section 6 in \cite{gut}), and because 
\[
     \sjj \Pi_{kj}\tq_j^Y(m,n)\Pi_{lj}\tq_l^X(m,n) \Eq M_X\umn(k,l) ,
\]
that
\beas
    \lefteqn{d_{TV}\Bigl(\law(X\umn(1) \giv X\umn(0)= e\uk),
             \otimes_{l=1}^K\Po(M_X\umn(k,l))\Bigr)}\\
     &&\Le \sjj\left(\frac{\Pi_{kj}}{\sqrt{mn}} 
       + \sqrt{\frac n m}\Pi_{kj}\tq_j^Y(m,n) \Blb \skk \tq_k^X(m,n)\Pi_{kj} \Brb^2
         \right),
\enas
an error of order $O(\sqrt{n/m})$.  Thus, in this regime, the
offspring distribution for the~$X\umn$ process, which determines the
distribution of~$W\umn$, approaches one with independent Poisson
components, having means given by the matrix~$M_X$. Again, this
entails the convergence of $U'{}\umn$ to~$U'$, but not the 
convergence of~$U\umn$.

Finally, if $m/n$ is small, the simple bound $(1-p)^l \ge 1-lp$ shows that 
\[
    \PP[Y_j\umn(1) \neq 0] \le (m/n)^{1/2}\Pi_{kj}\tq_j^Y(m,n),
\]
from which it follows that $\PP[W > 0] \le (m/n)^{1/2}\sjj \tq_j^Y(m,n)\Pi_{kj}$.
Hence, for $m/n \to 0$, the distance between two randomly chosen vertices
$(k_1,1)$ and $(k_2,1)$ is infinite, with probability close to~$1$.
However, if the two vertices $A$ and~$B$ do each have an edge joining them
to the object set, then each is connected to just one object with 
conditional probability of order $1 - O({\textstyle{\smn}})$, and 
the objects to which they are linked are distinct with probability
of order $1 - O(1/m)$. The distance between these two objects can now be
investigated, in this regime, by swapping the roles of vertices and objects,
and using the theorems above.

Thus if, in this scheme, $m/n$ converges to a finite or infinite 
limit, the
approximating probability distributions $\law(U'{\vphantom{X}}\umn)$ remain
relatively stable. In the error terms, the quantities $\t\umn$
and~$\g\umn$ converge to limits $\t$ and~$\g$, the corresponding 
quantities for the limit matrix~$M_X$.
The factor $n^{1/4}\erra^2$ behaves like $n^{-1/4}$, and $\tth_{i_0}$ like 
\hbox{$n^{-\d}\log n$,} for some $\d$ 
depending on~$M_X$, as long as $m/n$ is bounded below as~$\nti$;
in view of \Ref{gamma-eta} and~\Ref{theta-tilde-def}, it follows that
$\d \le 1/4$.
The discussion above shows that $m/n$  bounded below is the case of main interest.

\section{An exponential random graph model}\label{ER}
\setcounter{equation}{0}
Rank~$1$ matrices~$P = \a\b^T$ give rise to an exponential random graph model. 
The individual edges are independent, as before, and the probability of a 
vertex of type~$i$ 
connecting to an object of type~$j$ is of product form. These models 
have been extensively studied in the social science literature, see for example \cite{robinsetal}
and references therein; for applications to affiliation networks as bipartite networks, see for example \cite{wangetal}. 

In this case, 
\beas 
  M_X(k,l) &=& \sjj \a_k \b_j m_j \b_j \a_l n_l \Eq \a_k \a_l n_l \sjj m_j \b_j^2,
\enas  
so that 
\beas
    M_X \Eq C \a \a^T N_X,
\enas 
with $C = \sjj m_j \b_j^2$, has $\t = C\a^T N_X \a$, $\m = N_X\a/\bone^TN_X\a$ 
and $\n = C (\t^{-1}\bone^TN_X\a)\a$. Here, $\bone$ is a $K \times 1$-vector of 1's.

As shown in the preceding sections, the principal eigenvalue~$\t$ of~$M_X$ 
in our general multitype intersection graph is of 
critical importance in determining network distances.  It
turns out that its value 
can be bounded below by that obtained from an associated rank~$1$ 
matrix, adding to the importance of the exponential
models.  To see this, set $\BA_k := \sjj p_{kj}m_j$
to be the average degree of a type~$(k,1)$ 
vertex, and write $s_B^2 := \skk n_k (\BA_k)^2$.  

\begin{proposition}\label{e-value}
If the values $\BA_k$, $1\le k\le K$, are fixed, then 
\[
    \t\ \ge\ s_B^2 / m,
\]
and this value of~$\t$ is attained by taking $p_{kj} = \BA_k/m$
for all~$j$; with this choice of the $p_{kj}$'s, a vertex 
makes no distinction as to which types of object it has
links to.  The lower bound is minimized, if $\skk \BA_k = \BA_+$
is fixed, by taking $\BA_k = \BA_+/n$, so that all links have the
same probability~$p_{kj} = \BA_+/(mn)$.
\end{proposition}

\proof
The proof is taken from \cite{B78}, p.15.
The matrix $M_X := P N_Y P^T N_X$ has the same eigenvalues as the
symmetric matrix
\[
    V \ :=\ N_X^{1/2} P N_Y P^T N_X^{1/2}.
\]
Write $u_{kj} := (\BA_k)^{-1}mp_{kj} - 1$, and note that $\sjj m_j u_{kj}
= 0$ for each~$k$.  Then
\beas
   V_{ki} &=& \sqrt{n_k}\sjj p_{kj} m_j p_{ij} \sqrt{n_i} \\
   &=& m^{-1}\BA_k \sqrt{n_k} \BA_i \sqrt{n_i} \Blb
      1 + m^{-1}\sjj m_j u_{kj} u_{ij} \Brb.
\enas
From the Rayleigh-Ritz Theorem~\cite{HornJohnson}, Theorem~4.2.2, it follows that~$\t$, the largest eigenvalue of~$V$, is at least as large
as $e^T V e/(e^Te)$, for any~$e\in \re^K$. Taking $e_k
= \BA_k\sqrt{n_k}/s_B$, $1\le k\le K$, gives
\[
   \t \ \ge\ m^{-1} s_B^2 + \|v\|_2^2,
\]
where 
\[
   v_j \ :=\ m^{-1}\skk n_k(\BA_k)^2 u_{kj}\sqrt{m_j}, \qquad 1\le j\le J.
\]
Since, with $p_{kj} = \BA_k/m$ for all~$k$,  $V$ takes the form $m^{-1} w w^T$,
with $w = N_X^{1/2}\BA$, and hence has largest eigenvalue $m^{-1}\|w\|_2^2
= s_B^2/m$, this proves the first statement of the proposition; the
second is now immediate.
\ep

\nin Thus the value of~$\t$ for a given~$P$ is always bigger than 
that corresponding
to as homogeneous a choice of the link probabilities as is allowed
by the constraints on the average number of objects linked to a
given vertex. 

\medskip	 
Any rank one choice $P = \a\b^T$ has a minimality property, analogous
to that of Proposition~\ref{e-value}, but of a less intuitive
nature.  The matrix $P = \a\b^T$ minimizes the maximum 
eigenvalue of~$M_X$ among all choices of~$P$ satisfying the constraint
\[
    PN_Y\b \Eq (\b^T N_Y\b)\a.
\]

\vfil\eject

\end{document}